\documentclass[12pt]{article}
\usepackage{bbm}
\usepackage{amsfonts}
\usepackage{pifont}
\usepackage{hyperref}
\usepackage{mathrsfs}
\usepackage{indentfirst}
\usepackage{amsmath}
\usepackage{amssymb}
\usepackage[amsmath, thmmarks]{ntheorem}
\usepackage{cite}
\usepackage{graphicx}
\usepackage{tikz}
\newtheorem{definition}{\bf Definition}[section]
\newtheorem{lemma}{\bf Lemma}[section]
\newtheorem{theorem}{\bf Theorem}[section]
\newtheorem{remark}{\bf Remark}[section]
\newtheorem{corollary}{\bf Corollary}[section]
\newtheorem{proposition}{\bf Proposition}[section]
\newtheorem{example}{\bf Example}[section]

\def\proof{{\bf Proof.} }

\begin{document}
\setcounter{page}{1}

\title{
{\textbf{New structures for uninorms on bounded lattices}}\thanks {Supported by National Natural Science
Foundation of China (No.11871097)}}

\author{
Zhen-Yu Xiu$^1$\footnote{ E-mail address: \ xzy@cuit.edu.cn,\  xyz198202@163.com},  Yu-Xiu Jiang$^2$\footnote{Corresponding author. E-mail address:\ 191158845@qq.com }\\
\emph{\footnotesize  1,2  College of Applied Mathematics, Chengdu University of Information Technology,
}\\
\emph{
\footnotesize Chengdu 610000,  China }
}
\newcommand{\pp}[2]{\frac{ #1}{ #2}}
\date{}
\maketitle

\begin{quote}
{\bf Abstract}
In this article, we study new methods  for constructing uninorms on  bounded lattices.  First,  we present new methods  for constructing uninorms on  bounded lattices  under the additional constraints and prove that  some of these constraints are sufficient and necessary  for the uninorms.  Second,    we show  that   the additional constraints on  $t$-norms ($t$-conorms) and  $t$-subnorms ($t$-subconorms)  of  some  uninorms are exactly sufficient and necessary.  
 At last,  we give new constructions of uninorms  on arbitrary bounded lattices   by interation  based on $t$-conorms  and $t$-conorms.



\end{quote}

\section{\bf Introduction}

 The uninorms  on the unit interval $[0,1]$,  as a generalization of  $t$-norms and  $t$-conorms, are introduced by Yager and Rybalov \cite{RR96} and then applied to various fieleds \cite{BD99,MG09,MG11,WP07}, such as fuzzy logic, fuzzy set theory, expert systems, neural networks and so on.


 Recently,   because of  ``For general systems, where we cannot always expect real (or comparable) data, this extension of the underlying career together with its structure is significant..." \cite{GD21},  the researchers widely study uninorms on the bounded lattices instead of  the unit interval $[0,1]$.
About the methods for the construction of uninorms, they  mainly focus on t-norms (t-conorms)  \cite{EA21,SB14,GD19.2,GD19,GD20,GD17,GD16,YD19,FK15,AX20}, $t$-subnorms ( $t$-subconorms) \cite{XJ20,WJ21,HP21}, closure operators ( interior operators) \cite{GD21,XJ21,YO20,BZ21} and additive generators \cite{HeP}.

 For the extension,   as shown in many existing results, yielding uninorms on the bounded lattice  always needs  some additional constraints on the lattice,  $t$-norms (resp. $t$-conorms),  $t$-subnorms (resp. $t$-subconorms)  or  closure operators (resp. interior operators) to guarantee that the methods work well.   However, we may ask  the following question: are  these additional constraints necessary?

 As we see, a lot of constructions of uninorms were obtained by the existing methods. However, we may still  ask    are  there more methods  to construct uninorms on  bounded lattices,  different from the  existing  methods?
 In this paper, on one hand,  we  present  new methods  for constructing uninorms on the bounded lattices under some additional constraints and  show that  some constraints are necessary, that is,
 these constraints are sufficient and necessary for the uninorms.  Moreover,   some sufficient conditions of some  uninorms are modified to be  sufficient and necessary.
On the other hand,   we  present  new  methods  for constructing uninorms on arbitrary bounded lattices by the iterative method using $t$-norms and $t$-conorms, which are  different from the existing methods.  This method  is used to construct  uninorms on arbitrary bounded lattices for the first time.

%


%
The rest of this paper is organized as follows. In section $2$, we recall some preliminaries. In section $3$, we construct uninorms based on the sufficient and necessary  conditions of the $t$-subnorm ($t$-subconorm) and  the sufficient and necessary constraints of bounded lattices, respectively.  In section $4$, we give the modifications of some uninorms.
In Section $5$,  we  present  a new method  for constructing uninorms on arbitrary bounded lattices by the iterative method using $t$-norms and $t$-conorms.
In Section $6$, some conclusions are made.

\section{\bf Preliminaries }
\vspace{.2 cm}
In this section, we recall some preliminary details concerning bounded lattice and results related to aggregation functions on bounded lattices.
\begin{definition}[\cite{GB67}]\label{de2.1}
\emph{A lattice $(L,\leq)$  is bounded if it has top and bottom elements, which are written as $1$  and $0$, respectively, that is, there exit two elements $1,0\in L$  such that $0\leq x\leq 1$ for all $x\in L$.}\\

\emph{Throughout this article, unless stated otherwise, we denote $L$  as a bounded lattice with the top and bottom elements  $1$ and $0$, respectively.}
\end{definition}

\begin{definition}[\cite{GB67}]\label{de2.2}
\emph{Let $L$ be a bounded lattice, $a,b\in L$ with $a\leq b$. A subinterval $[a,b]$ of $L$ is defined as
$$[a, b]=\{x\in L: a\leq x \leq b\}.$$
Similarly, we can define $[a, b[=\{x\in L: a\leq x < b\}, ]a, b]=\{x\in L: a< x \leq b\}$ and $]a, b[=\{x\in L: a< x < b\}$. Denote $D(e)=I_{e}\times \{1\}\cup \{1\}\times I_{e}\cup [0,e]\times \{1\}\cup \{1\}\times [0,e]\cup [e,1]\times \{e\}\cup \{e\}\times [e,1]$, $D(e)'=I_{e}\times \{1\}\cup \{1\}\times I_{e}\cup ]0,e[\times \{1\}\cup \{1\}\times ]0,e[\cup ]e,1]\times \{e\}\cup \{e\}\times ]e,1]$, $E(e)=I_{e}\times \{0\}\cup \{0\}\times I_{e}\cup [e,1]\times \{0\}\cup \{0\}\times [e,1]\cup [0,e]\times \{e\}\cup \{e\}\times [0,e]$, $E(e)'=I_{e}\times \{0\}\cup \{0\}\times I_{e}\cup ]e,1[\times \{0\}\cup \{0\}\times ]e,1[\cup [0,e[\times \{e\}\cup \{e\}\times [0,e[$. If $a$  and $b$  are incomparable, then we use the notation $a\parallel b$. For $e\in L\setminus\{0,1\}$, we denote the set of all incomparable elements with $e$  by $I_{e}$, that is, $I_{e}=\{x\in L \mid x \parallel e\}$.}
\end{definition}

\begin{definition}[\cite{EP04}]\label{de2.3}
\emph{Let $(L,\leq,0,1)$  be a bounded lattice.}\\\
\emph{(i) An operation $T:L^{2}\rightarrow L$  is called a $t$-norm on $L$ if it is commutative, associative, and increasing with respect to both variables, and it has the neutral element $1\in L$, that is, $T(1,x)=x$  for all $x\in L$. The greatest $t$-norm $T_{\wedge}$  on $L$  is given as $T_{\wedge}(x,y)=x\wedge y$;}\\\
\emph{(ii) An operation $S:L^{2}\rightarrow L$  is called a $t$-conorm on  $L$ if it is commutative, associative, and increasing with respect to both variables, and it has the neutral element $0\in L$, that is, $S(0,x)=x$  for all $x\in L$. The smallest t-conorm $S_{\vee}$  on   $L$ is given as $S_{\vee}(x,y)=x\vee y$.}\
\end{definition}

\begin{definition}[\cite{FK15}]\label{de2.4}
\emph{Let $(L,\leq,0,1)$  be a bounded lattice. An operation $U:L^{2}\rightarrow L$  is called a uninorm on  $L$ (a uninorm if  $L$ is fixed) if it is commutative, associative, and increasing with respect to both variables, and it has the neutral element $e\in L$, that is, $U(e,x)=x$  for all $x\in L$.}
\end{definition}


\begin{proposition}[\cite{FK15}]\label{pro2.1}\rm
Let $(L,\leq,0,1)$  be a bounded lattice and $U$ be a uninorm on $L$ with a neutral element $e\in L\setminus\{0,1\}$. Then we have the following:\\
(i)  $T_{e}=U\mid [0,e]^2:[0,e]^2\rightarrow [0,e] $ is a $t$-norm on $[0,e]$;\\
(ii)  $S_{e}=U\mid [e,1]^2:[e,1]^2\rightarrow [e,1]$ is a $t$-conorm on $[e,1]$.\\

\emph{$T_{e}$ and $S_{e}$  given in proposition \ref{pro2.1}  are called the underlying $t$-norm and $t$-conorm of a uninorm $U$  on a bounded lattice $L$  with the neutral element $e$, respectively. Throughout this study, we denote $T_{e}$  as the underlying $t$-norm and $S_{e}$  as the underlying $t$-conorm of a given uninorm  $U$ on $L$.}
\end{proposition}

\begin{proposition}[\cite{WJ21}]\label{pro2.2}
\emph{Let $S$ be a nonempty set and $A_{1},A_{2},\ldots,A_{n}$ be subsets of $S$. Let $H$ be a commutative binary operation on $S$. Then $H$ is associative on $\bigcup \limits_{i\in I}A_{i}$ if the following statements hold:}\\\
\emph{(i) $H$ is alternating associative on $(A_{i},A_{j},A_{k})$ for every combination $\{i,j,k\}$ of $size $ $3$ chosen from $\{1,2,\ldots,n\}$;}\\\
\emph{(ii) $H$ is alternating associative on $(A_{i},A_{i},A_{j})$ for every combination $\{i,j\}$ of $size $ $2$ chosen from $\{1,2,\ldots,n\}$;}\\\
\emph{(iii) $H$ is alternating associative on $(A_{i},A_{j},A_{j})$ for every combination $\{i,j,j\}$ of $size$ $2$ chosen from $\{1,2,\ldots,n\}$;}\\\
\emph{(iv) $H$ is alternating associative on $(A_{i},A_{i},A_{i})$ for every $i \in \{1,2,\ldots,n\}$.}
\end{proposition}

\begin{definition}[\cite{EP04}]\label{de2.4}
\emph{Let $(L,\leq ,0,1)$  be a bounded lattice.}\\\
 \emph{(i) An operation $F:L^{2}\rightarrow L$  is called a $t$-subnorm on $L$  if it is commutative, associative, increasing with respect to both variables, and $F(x,y)\leq x\wedge y$  for all $x,y\in L$.}\\\
\emph{(ii) An operation $R:L^{2}\rightarrow L$  is called a $t$-subconorm on $L$  if it is commutative, associative, increasing with respect to both variables, and $R(x,y)\geq x\vee y$  for all $x,y\in L$.}\
\end{definition}

\begin{lemma}[\cite{GD19}]\label{le2.1}
Let $L$ be a bounded lattice with $e\in L\setminus\{0,1\}$. \\
(1)  If $x\parallel y$  for all $x\in I_{e}$  and $y\in[e,1[$, then $z\vee t\in I_e$ or $z\vee t=1$ for all $z,t\in I_e$.\\
(2) If $x\parallel y$  for all $x\in I_{e}$  and $y\in]0,e]$, then $z\wedge t\in I_e$ or $z\wedge t=0$ for all $z,t\in I_e$.
\end{lemma}

\begin{remark}\label{Re2.8}  By lemma \ref{le2.1}, we can get the following results. \\
(1)  If $x\parallel y$  for all $x\in I_{e}$  and $y\in [e,1[ $, then $x\vee y\in I_{e}$  or $x\vee y=1$ for all ${x,y}\in I_{e}$. \\
(2)  If $x\parallel y$  for all $x\in I_{e}$  and $y\in ]0,e] $, then $x\wedge y\in I_{e}$  or $x\wedge y=0$ for all ${x,y}\in I_{e}$. \\
(3)   If $x\parallel y$  for all $x\in I_{e}$  and $y\in [e,1[ $, then $R(x,y)\in I_{e}$  or $R(x,y)=1$ for all ${x,y}\in I_{e}$.\\
(4)  If $x\parallel y$  for all $x\in I_{e}$  and $y\in ]0,e] $, then $F(x,y)\in I_{e}$  or $F(x,y)=0$ for all ${x,y}\in I_{e}$.\\

\end{remark}

\begin{theorem}[\cite{YD20}]\label{th2.0}
\emph{Let $e\in L\setminus\{0,1\}$, $T_{e}$ is a t-norm on $[0,e]$ and $S_{e}$ is a t-norm on $[e,1]$. Then the function $U_{(T,e)}:L^{2}\rightarrow L$ and $U_{(S,e)}:L^{2}\rightarrow L$ defined as follows
}

\emph{$U_{(T,e)}(x,y)=\begin{cases}
T_{e}(x, y) &\mbox{if } (x,y)\in [0,e]^{2},\\
x &\mbox{if } (x,y)\in I_{e}\times [0,e],\\
y &\mbox{if } (x,y)\in [0,e]\times I_{e},\\
x\vee y\vee e &\mbox{}otherwise,\\
\end{cases}$\\
and}

\emph{$U_{(S,e)}(x,y)=\begin{cases}
S_{e}(x, y) &\mbox{if } (x,y)\in [e,1]^{2},\\
x &\mbox{if } (x,y)\in I_{e}\times [e,1],\\
y &\mbox{if } (x,y)\in [e,1]\times I_{e},\\
x\wedge y\wedge e  &\mbox{}otherwise.\\
\end{cases}$}
\end{theorem}

\begin{theorem}[\cite{GD19}]\label{th2.1}
Let $(L,\leq,0,1)$  be a bounded lattice with $e\in L\setminus\{0,1\}$.\\
(i) Let $T_{e}$  be a $t$-norm on $[0,e]$ such that $T_{e}(x,y)>0$ for all $x,y>0$ and $S_{e}$  be a $t$-conorm on $[e,1]$ such that $S_{e}(x,y)<1$ for all $x,y<1$. If $y\parallel x$  for all $x\in I_{e}$ and $y\in[e,1[$, then the function $U_{1}^{e}:L^{2}\rightarrow L$  is a uninorm on $L$  with the neutral element $e$, where

$U_{1}^{e}(x,y)=\begin{cases}
T_{e}(x, y) &\mbox{if } (x,y)\in [0,e]^{2},\\
S_{e}(x, y) &\mbox{if } (x,y)\in [e,1]^{2},\\
x &\mbox{if } (x,y)\in I_{e}\times [e,1[\cup I_{e}\times ]0,e[,\\
y &\mbox{if } (x,y)\in [e,1[\times I_{e}\cup ]0,e[\times I_{e},\\
x\vee y &\mbox{if } (x,y)\in I_{e}^2\cup I_{e}\times \{1\}\cup \{1\}\times I_{e}\cup ]0,e[\times \{1\}\cup \{1\}\times ]0,e[,\\
x\wedge y &\mbox{  } otherwise.\\
\end{cases}$\\
(ii) Let $T_{e}$  be a $t$-norm on $[0,e]$ such that $T_{e}(x,y)>0$ for all $x,y>0$ and $S_{e}$  be a $t$-conorm on $[e,1]$ such that $S_{e}(x,y)<1$ for all $x,y<1$. If $y\parallel x$  for all $x\in I_{e}$ and $y\in]0,e]$, then the function $U_{2}^{e}:L^{2}\rightarrow L$  is a uninorm on $L$  with the neutral element $e$, where

$U_{2}^{e}(x,y)=\begin{cases}
T_{e}(x, y) &\mbox{if } (x,y)\in [0,e]^{2},\\
S_{e}(x, y) &\mbox{if } (x,y)\in [e,1]^{2},\\
x &\mbox{if } (x,y)\in I_{e}\times ]e,1[\cup I_{e}\times ]0,e],\\
y &\mbox{if } (x,y)\in ]e,1[\times I_{e}\cup ]0,e]\times I_{e},\\
x\wedge y &\mbox{if } (x,y)\in I_{e}^2\cup I_{e}\times \{0\}\cup \{0\}\times I_{e}\cup ]e,1[\times \{0\}\cup \{0\}\times ]e,1[,\\
x\vee y &\mbox{  } otherwise.\\
\end{cases}$\\
\end{theorem}

\begin{theorem}[\cite{XJ20}]\label{th2.2}
Let $(L,\leq,0,1)$  be a bounded lattice with $e\in L\setminus\{0,1\}$.\\
(i) Let $T_{e}$  be a $t$-norm on $[0,e]$ and $R$ is a $t$-subconorm on $L$ such that $R(x,y)<1$ for all $x,y<1$. If $x\parallel y$  for all $x\in I_{e}$  and $y\in [e,1[$,  then the function $U_{R}^{e}:L^{2}\rightarrow L$  is a uninorm on $L$  with the neutral element $e$, where

$U_{R}^{e}(x,y)=\begin{cases}
T_{e}(x, y) &\mbox{if } (x,y)\in [0,e]^{2},\\
x &\mbox{if } (x,y)\in I_{e}\times [0,e]\cup I_{e}\times [e,1[,\\
y &\mbox{if } (x,y)\in [0,e]\times I_{e}\cup [e,1[\times I_{e},\\
R(x, y) &\mbox{if } (x,y)\in I_{e}^2\cup ]e,1]^{2},\\
x\vee y &\mbox{if } (x,y)\in D(e),\\
x\wedge y &\mbox{  } otherwise.\\
\end{cases}$\\
(ii) Let  $S_{e}$ be a $t$-conorm on $[e,1]$ and $F$ is a $t$-subnorm on $L$ such that $F(x,y)>0$ for all $x,y>0$. If $x\parallel y$  for all $x\in I_{e}$  and $y\in ]0,e]$, then the function $U_{F}^{e}:L^{2}\rightarrow L$  is a uninorm on $L$  with the neutral element $e$, where

$U_{F}^{e}(x,y)=\begin{cases}
S_{e}(x, y) &\mbox{if } (x,y)\in [e,1]^{2},\\
x &\mbox{if } (x,y)\in I_{e}\times ]0,e]\cup I_{e}\times [e,1],\\
y &\mbox{if } (x,y)\in ]0,e]\times I_{e}\cup [e,1]\times I_{e},\\
F(x, y) &\mbox{if } (x,y)\in I_{e}^2\cup [0,e[^{2},\\
x\wedge y &\mbox{if } (x,y)\in E(e),\\
x\vee y &\mbox{  } otherwise.\\
\end{cases}$\\
\begin{theorem}[\cite{HP21}]\label{th2.3}
Let $(L,\leq,0,1)$  be a bounded lattice with $e\in L\setminus\{0,1\}$.\\
(i) Let $T_{e}$  be a $t$-norm on $[0,e]$ and $R$  be a $t$-subconorm on $L\setminus[0,e]$, then the function $U_{R}:L^{2}\rightarrow L$ defined by

$U_{R}(x,y)=\begin{cases}
T_{e}(x, y) &\mbox{if } (x,y)\in [0,e]^{2},\\
x &\mbox{if } (x,y)\in [0,e[\times (L\setminus[0,e]),\\
y &\mbox{if } (x,y)\in (L\setminus[0,e])\times [0,e[,\\
y &\mbox{if } (x,y)\in \{e\}\times (L\setminus[0,e]),\\
x &\mbox{if } (x,y)\in (L\setminus[0,e])\times \{e\},\\
R(x,y) &\mbox{if } (x,y)\in (L\setminus[0,e])^{2},
\end{cases}$\\
is a uninorm on $L$  with $e\in L\setminus\{0,1\}$ iff $y< x$  for all $x\in I_{e}$ and $y\in]0,e[$.\\
(ii) Let $S_{e}$  be a $t$-conorm on $[e,1]$ and $F$  be a $t$-subnorm on $L\setminus[e,1]$, then the function $U_{F}:L^{2}\rightarrow L$ defined by

 $U_{F}(x,y)=\begin{cases}
S_{e}(x, y) &\mbox{if } (x,y)\in [e,1]^{2},\\
x &\mbox{if } (x,y)\in ]e,1]\times (L\setminus[e,1]),\\
y &\mbox{if } (x,y)\in (L\setminus[e,1])\times ]e,1],\\
y &\mbox{if } (x,y)\in \{e\}\times (L\setminus[e,1]),\\
x &\mbox{if } (x,y)\in (L\setminus[e,1])\times \{e\},\\
F(x,y) &\mbox{if } (x,y)\in (L\setminus[e,1])^{2},
\end{cases}$\\
is a uninorm on $L$  with $e\in L\setminus\{0,1\}$ iff $x< y$  for all $x\in I_{e}$ and $y\in ]e,1[$.

\end{theorem}

\begin{theorem}[\cite{GD21}]\label{th2.4}
Let $(L,\leq,0,1)$  be a bounded lattice with $e\in L\setminus\{0,1\}$.\\
(i) Let $T_{e}$  be a $t$-norm on $[0,e]$, then the function $U_{\wedge}:L^{2}\rightarrow L$ defined by

$U_{\wedge}(x,y)=\begin{cases}
T_{e}(x, y) &\mbox{if } (x,y)\in [0,e]^{2},\\
x\wedge y &\mbox{if } (x,y)\in [0,e[\times I_{e}\cup I_{e}\times [0,e[ \cup [0,e[\times [e,1]\cup [e,1]\times [0,e[,\\
y &\mbox{if } (x,y)\in \{e\}\times I_{e},\\
x &\mbox{if } (x,y)\in I_{e}\times \{e\},\\
x\vee y & otherwise,
\end{cases}$\\
is a uninorm on $L$  with $e\in L\setminus\{0,1\}$ iff $x> y$  for all $x\in I_{e}$  and $y\in [0,e[$.\\
(ii) Let $S_{e}$ be a $t$-conorm on $[e,1]$,  then the function $U_{\vee}:L^{2}\rightarrow L$ defined by

$U_{\vee}(x,y)=\begin{cases}
S_{e}(x, y) &\mbox{if } (x,y)\in [e,1]^{2},\\
x\vee y &\mbox{if } (x,y)\in ]e,1]\times I_{e}\cup I_{e}\times ]e,1] \cup [0,e[\times ]e,1]\cup ]e,1]\times [0,e[,\\
y &\mbox{if } (x,y)\in \{e\}\times I_{e},\\
x &\mbox{if } (x,y)\in I_{e}\times \{e\},\\
x\wedge y & otherwise,
\end{cases}$\\
is a uninorm on $L$  with $e\in L\setminus\{0,1\}$ iff $x< y$  for all $x\in I_{e}$ and $y\in ]e,1]$.\\

\end{theorem}
\end{theorem}

\begin{theorem}[\cite{YD20}]\label{th2.5}
\emph{Assume that $T: L^{2}\rightarrow L$ and $S: L^{2}\rightarrow L$ are a t-norm and a t-conorm, respectively. Let $e\in L$ be an element distinct from both $0_{L}$ and $1_{L}$. Let us denote by $U_{d}: L^{2}\rightarrow L$ and $U_{c}: L^{2}\rightarrow L$ the following:}

\emph{$U_{d}(x,y)=\begin{cases}
T(x, y) &\mbox{if } (x,y)\in [0,e]^{2},\\
S(x, y) &\mbox{if } (x,y)\in (e,1]^{2},\\
x &\mbox{if } (x,y)\in (I_{e}\cup(e,1])\times [0,e],\\
y &\mbox{if } (x,y)\in [0,e]\times (I_{e}\cup(e,1]),\\
S(x\vee a, y\vee a) &\mbox{}otherwise,\\
\end{cases}$\\
and}

\emph{$U_{c}(x,y)=\begin{cases}
T(x, y) &\mbox{if } (x,y)\in [0,e)^{2},\\
S(x,y) &\mbox{if } (x,y)\in [e,1]^{2},\\
x &\mbox{if } (x,y)\in (I_{e}\cup [0,e))\times [e,1],\\
y &\mbox{if } (x,y)\in [e,1]\times (I_{e}\cup [0,e)),\\
T(x\wedge a, y\wedge a) &\mbox{}otherwise.\\
\end{cases}$\\}
\end{theorem}

\section{ New construction of uninorms on bounded lattices}

  In this section,  we present  new   construction of uninorms  on some appropriate bounded lattices $L$, which are constructed by means of $t$-conorms ($t$-norms) and  $t$-subconorms ($t$-subnorms).  Moreover, the sufficient and necessary conditions are provided for these uninorms.

\begin{theorem}\label{th3.1}
Let $(L,\leq,0,1)$ be a bounded lattice with $e \in L\setminus \{0,1\}$, $T_{e}$ be a $t$-norm on $[0,e]$ and $S_{e}$ be a $t$-conorm on $[e,1]$.\\
$(i)$ If $x\parallel y$  for all $x\in I_{e}$ and $y\in [e,1[$, then the function $U_1:L^{2}\rightarrow L$ defined by

$U_{1}(x,y)=\begin{cases}
T_{e}(x, y) &\mbox{if } (x,y)\in [0,e]^{2},\\
S_{e}(x, y) &\mbox{if } (x,y)\in [e,1]^{2},\\
x &\mbox{if } (x,y)\in I_{e}\times [e,1[\cup I_{e}\times [0,e[,\\
y &\mbox{if } (x,y)\in [e,1[\times I_{e}\cup [0,e[\times I_{e},\\
x\vee y &\mbox{if } (x,y)\in {I_{e}}^{2}\cup I_{e}\times \{1\}\cup \{1\}\times I_{e}\cup [0,e[\times \{1\}\cup \{1\}\times [0,e[,\\
x\wedge y &\mbox{} otherwise,
\end{cases}$\\
is a uninorm on $L$  with $e\in L\setminus\{0,1\}$ iff $S_{e}(x,y)<1$ for all $x,y<1$.\\
$(ii)$ If $x\parallel y$  for all $x\in I_{e}$ and $y\in ]0,e]$, then the function $U_2:L^{2}\rightarrow L$ defined by

$U_{2}(x,y)=\begin{cases}
T_{e}(x, y) &\mbox{if } (x,y)\in [0,e]^{2},\\
S_{e}(x, y) &\mbox{if } (x,y)\in [e,1]^{2},\\
x &\mbox{if } (x,y)\in I_{e}\times ]e,1]\cup I_{e}\times ]0,e],\\
y &\mbox{if } (x,y)\in ]e,1]\times I_{e}\cup ]0,e]\times I_{e},\\
x\wedge y &\mbox{if } (x,y)\in {I_{e}}^{2}\cup I_{e}\times \{0\}\cup \{0\}\times I_{e}\cup ]e,1]\times \{0\}\cup \{0\}\times ]e,1],\\
x\vee y &\mbox{} otherwise,
\end{cases}$\\
is a uninorm on $L$  with $e\in L\setminus\{0,1\}$ iff $T_{e}(x,y)>0$ for all $x,y>0$.

\end{theorem}
\proof
Let $(L,\leq,0,1)$ be a bounded lattice with $e \in L\setminus \{0,1\}$. If $x\parallel y$  for all $x\in I_{e}$ and $y\in [e,1[$, then we give the proof of the fact that $U_{1}$ is a uninorm iff $S_{e}(x,y)<1$ for all $x,y<1$. The same result for $U_{2}$ can be obtained by similar arguments.

Necessity.
Let   $U_{1}$ be a uninorm on $L$ with the identity $e$. Then  we prove that $S_{e}(x,y)<1$ for all $x,y<1$. Assume that there are some elements $x\in ]e,1[$ and $y\in ]e,1[$ such that $S_{e}(x,y)=1$. If $z\in I_{e}$, then  $U_{1}(x,U_{1}(y, z))=U_{1}(x,z)=z$ and $U_{1}(U_{1}(x,y), z)=U_{1}(S_{e}(x,y),z)=U_{1}(1,z)=1$. Since $S_{e}(x,y)=1$, the associativity property of  $U_{1}$ is violated. Then $U_{1}$ is not a uninorm on $L$, which is a contradiction. Hence, $S_{e}(x,y)<1$ for all $x,y<1$.

Sufficiency.
To prove the sufficiency, we compare $U_{1}$ with $U_{1}^{e}$ in  Theorem $\ref{th2.1}$  and find that
besides the condition  that $T_{e}(x,y)>0$ for all $x,y>0$,  the values of $U_{1}$ and $U_{1}^{e}$  are different on $I_{e} \times \{0\}\cup \{1\}\times \{0\}$, $\{0\}\times I_{e}\cup \{0\} \times\{1\}$ and $]0,e[\times ]0,e[$. So, combined with  Theorem $\ref{th2.1}$, the proof of the sufficiency  can be reduced to the above areas.
The commutativity of $U_{1}$  and the fact that $e$  is a neutral element are evident. Hence, we prove only the monotonicity and the associativity of $U_{1}$.

I. Monotonicity: We prove that if $x\leq y$, then for all $z\in L$, $U_{1}(x,z)\leq U_{1}(y,z)$.

If one of $x$, $y$, $z$ is equal to $1$, then it is clear that $U_{1}(x,z)\leq U_{1}(y,z)$. There is no $1$ in any of the following cases.

1. Let $x\in [0,e[$.

\ \ \ 1.1. $y\in [0,e[$,

\ \ \ \ \ \ 1.1.1. $z\in I_{e}$,

\ \ \ \ \ \ \ \ \ \ \ \ $U_1(x,z)=z=U_1(y,z)$

\ \ \ 1.2. $y\in [e,1[$,

\ \ \ \ \ \ 1.2.1. $z\in I_{e}$,

\ \ \ \ \ \ \ \ \ \ \ \ $U_1(x,z)=z=U_1(y,z)$

\ \ \ 1.3. $y\in I_{e}$,

\ \ \ \ \ \ 1.3.1. $z\in [0,e[$,

\ \ \ \ \ \ \ \ \ \ \ \ $U_1(x,z)=T_{e}(x,z)\leq y=U_1(y,z)$

\ \ \ \ \ \ 1.3.2. $z\in [e,1[$,

\ \ \ \ \ \ \ \ \ \ \ \ $U_1(x,z)=x\leq y=U_1(y,z)$

\ \ \ \ \ \ 1.3.3. $z\in I_{e}$,

\ \ \ \ \ \ \ \ \ \ \ \ $U_1(x,z)=z\leq y\vee z=U_1(y,z)$

2. Let $x\in I_{e}$.

\ \ \ 2.1. $y\in I_{e}$,

\ \ \ \ \ \ 2.1.1. $z\in [0,e[$,

\ \ \ \ \ \ \ \ \ \ \ \ $U_1(x,z)=x\leq y=U_1(y,z)$

II. Associativity: We need to prove that $U_1(x,U_1(y,z))=U_1(U_1(x,y),z)$  for all $x,y,z\in L$.

If one of $x,y,z$ is equal to $1$, then it is clear that $U_1(x,U_1(y,z))=1=U_1(U_1(x,y),z)$. So there is no $1$ in any of the following cases.

1. Let $x\in [0,e[$.

\ \ \ 1.1. $y\in [0,e[$,

\ \ \ \ \ \ 1.1.1. $z\in I_{e}$,

\ \ \ \ \ \ \ \ \ \ \ \ $U_1(x,U_1(y,z))=U_1(x,z)=z=U_1(T_{e}(x,y),z)=U_1(U_1(x,y),z)$

\ \ \ 1.2. $y\in [e,1[$,

\ \ \ \ \ \ 1.2.1. $z\in I_{e}$,

\ \ \ \ \ \ \ \ \ \ \ \ $U_1(x,U_1(y,z))=U_1(x,z)=z=U_1(x,z)U_1(U_1(x,y),z)$

\ \ \ 1.3. $y\in I_{e}$,

\ \ \ \ \ \ 1.3.1. $z\in [0,e[$,

\ \ \ \ \ \ \ \ \ \ \ \ $U_1(x,U_1(y,z))=U_1(x,y)=y=U_1(y,z)=U_1(U_1(x,y),z)$

\ \ \ \ \ \ 1.3.2. $z\in [e,1[$,

\ \ \ \ \ \ \ \ \ \ \ \ $U_1(x,U_1(y,z))=U_1(x,y)=y=U_1(y,z)=U_1(U_1(x,y),z)$

\ \ \ \ \ \ 1.3.3. $z\in I_{e}$,

\ \ \ \ \ \ \ \ \ 1.3.3.1. $y\vee z\in I_{e}$,

\ \ \ \ \ \ \ \ \ \ \ \ $U_1(x,U_1(y,z))=U_1(x,y\vee z)=y\vee z=U_1(y,z)=U_1(U_1(x,y),z)$

\ \ \ \ \ \ \ \ \ 1.3.3.2. $y\vee z=1$,

\ \ \ \ \ \ \ \ \ \ \ \ $U_1(x,U_1(y,z))=U_1(x,y\vee z)=1=y\vee z=U_1(y,z)=U_1(U_1(x,y),z)$

2. Let $x\in [e,1[$.

\ \ \ 2.1. $y\in [0,e[$,

\ \ \ \ \ \ 2.1.1. $z\in I_{e}$,

\ \ \ \ \ \ \ \ \ \ \ \ $U_1(x,U_1(y,z))=U_1(x,z)=z=U_1(y,z)=U_1(U_1(x,y),z)$

\ \ \ 2.2. $y\in I_{e}$,

\ \ \ \ \ \ 2.2.1. $z\in [0,e[$,

\ \ \ \ \ \ \ \ \ \ \ \ $U_1(x,U_1(y,z))=U_1(x,y)=y=U_1(y,z)=U_1(U_1(x,y),z)$

3. Let $x\in I_{e}$.

\ \ \ 3.1. $y\in [0,e[$,

\ \ \ \ \ \ 3.1.1. $z\in [0,e[$,

\ \ \ \ \ \ \ \ \ \ \ \ $U_1(x,U_1(y,z))=U_1(x,T_{e}(y,z))=x=U_1(x,z)=U_1(U_1(x,y),z)$

\ \ \ \ \ \ 3.1.2. $z\in [e,1[$,

\ \ \ \ \ \ \ \ \ \ \ \ $U_1(x,U_1(y,z))=U_1(x,y)=x=U_1(x,z)=U_1(U_1(x,y),z)$

\ \ \ \ \ \ 3.1.3. $z\in I_{e}$,

\ \ \ \ \ \ \ \ \ \ \ \ $U_1(x,U_1(y,z))=U_1(x,z)=x\vee z=U_1(x,z)=U_1(U_1(x,y),z)$

\ \ \ 3.2. $y\in [e,1[$,

\ \ \ \ \ \ 3.2.1. $z\in [0,e[$,

\ \ \ \ \ \ \ \ \ \ \ \ $U_1(x,U_1(y,z))=U_1(x,z)=x=U_1(x,z)=U_1(U_1(x,y),z)$

\ \ \ 3.3. $y\in I_{e}$,

\ \ \ \ \ \ 3.3.1. $z\in [0,e[$,

\ \ \ \ \ \ \ \ \ 3.3.1.1. $x\vee y\in I_{e}$,

\ \ \ \ \ \ \ \ \ \ \ \ $U_1(x,U_1(y,z))=U_1(x,y)=x\vee y=U_1(x\vee y,z)=U_1(U_1(x,y),z)$

\ \ \ \ \ \ \ \ \ 3.3.1.2. $x\vee y=1$,

\ \ \ \ \ \ \ \ \ \ \ \ $U_1(x,U_1(y,z))=U_1(x,y)=x\vee y=1=U_1(x\vee y,z)=U_1(U_1(x,y),z)$

Hence, $U_{1}$ is a uninorm on $L$ with a neutral element $e$.\\

It is worth pointing out that if the condition that $x\parallel y$  for all $x\in I_{e}$  and $y\in [e,1[$  in Theorem \ref{th3.1} is not satisfied, then $U_{1}$  may not be a uninorm on $L_{1}$. Similarly, if the condition that $x\parallel y$  for all $x\in I_{e}$  and $y\in ]0,e]$  in Theorem \ref{th3.1} is not satisfied, then $U_{2}$  may not be a uninorm on $L_{1}$ . The following example will show the above facts.

\begin{example}\label{ex3.1}
 Consider the bounded lattice $L_{1}=\{0,a,e,c,f,g,1\}$   depicted by the Hasse diagram in Fig.$1$, where $a\in [0,e[$, $g\in ]e,1]$, $c\in I_{e}$, $f\in I_{e}$, $g>c$, $g>f$,$a<c$ and $a<f$.
Take $T_{e}(x,y)=T_{\wedge}(x,y)$  on  $[0,e]$ and $S_e(x,y)=S_{\vee}(x,y)$ on $[e,1]$. The function $U_{1}$  on $L_{1}$, shown in Table \ref{Tab:01}, is not a uninorm on $L_{1}$  with the neutral element $e$. In fact, the associativity is not satisfied, since $U_{1}(a,U_{1}(c,f))=U_{1}(a,g)=a$ and $U_{1}(U_{1}(a,c),f)=U_{1}(c,f)=g$. Similarly, the function $U_{2}$  on $L_{1}$, shown in Table \ref{Tab:02}, is not a uninorm on $L_{1}$  with the neutral element $e$. In fact, the associativity is not satisfied, since $U_{2}(g,U_{2}(c,f))=U_{2}(g,a)=g$ and $U_{2}(U_{2}(g,c),f)=U_{2}(c,f)=a$.

\par\noindent\vskip50pt

\begin{minipage}{11pc}
\setlength{\unitlength}{0.75pt}
\begin{picture}(600,180)
\put(240,40){\circle{5}}\put(238,29){\makebox(0,0)[l]{\footnotesize$0$}}
\put(240,80){\circle{5}}\put(245,80){\makebox(0,0)[l]{\footnotesize$a$}}
\put(200,120){\circle{5}}\put(189,120){\makebox(0,0)[l]{\footnotesize$c$}}
\put(240,120){\circle{5}}\put(245,120){\makebox(0,0)[l]{\footnotesize$e$}}
\put(280,120){\circle{5}}\put(285,120){\makebox(0,0)[l]{\footnotesize$f$}}
\put(240,160){\circle{5}}\put(245,164){\makebox(0,0)[l]{\footnotesize$g$}}
\put(240,200){\circle{5}}\put(237,213){\makebox(0,0)[l]{\footnotesize$1$}}
\put(240,43){\line(0,1){34}}
\put(240,83){\line(0,1){34}}
\put(240,123){\line(0,1){34}}
\put(240,163){\line(0,1){34}}
\put(238,82){\line(-1,1){36}}
\put(242,82){\line(1,1){36}}
\put(202,122){\line(1,1){36}}
\put(278,122){\line(-1,1){36}}
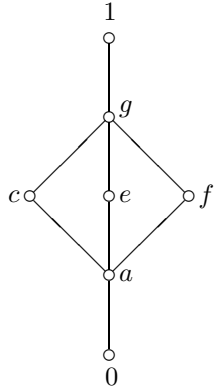
\put(170,-10){\emph{Fig.1. The lattice $L_{1}$}}
\end{picture}
\end{minipage}

\begin{table}[htpp]
\centering
\caption{The function $U_{1}$  on $L_{1}$  given in Fig.1.}
\ \ \ \ \\
\label{Tab:01}
\begin{tabular}{c|ccccccc}

 $U_{1}$ & $0$ & $a$ & $e$ & $c$ & $f$ & $g$ & $1$\\
 \hline
$0$ & $0$ & $0$ & $0$ & $c$ & $f$ & $0$ & $1$\\
$a$ & $0$ & $a$ & $a$ & $c$ & $f$ & $a$ & $1$\\
$e$ & $0$ & $a$ & $e$ & $c$ & $f$ & $g$ & $1$\\
$c$ & $c$ & $c$ & $c$ & $c$ & $g$ & $c$ & $1$\\
$f$ & $f$ & $f$ & $f$ & $g$ & $f$ & $f$ & $1$\\
$g$ & $0$ & $a$ & $g$ & $c$ & $f$ & $g$ & $1$\\
$1$ & $1$ & $1$ & $1$ & $1$ & $1$ & $1$ & $1$\\
\end{tabular}
\end{table}

\par\noindent\vskip50pt

\begin{table}[htpp]
\centering
\caption{The function $U_{2}$  on $L_{1}$  given in Fig.1.}
\ \ \ \ \\
\label{Tab:02}
\begin{tabular}{c|ccccccc}

 $U_{2}$ & $0$ & $a$ & $e$ & $c$ & $f$ & $g$ & $1$\\
 \hline
$0$ & $0$ & $0$ & $0$ & $0$ & $0$ & $0$ & $0$\\
$a$ & $0$ & $a$ & $a$ & $c$ & $f$ & $g$ & $1$\\
$e$ & $0$ & $a$ & $e$ & $c$ & $f$ & $g$ & $1$\\
$c$ & $0$ & $c$ & $c$ & $c$ & $a$ & $c$ & $c$\\
$f$ & $0$ & $f$ & $f$ & $a$ & $f$ & $f$ & $f$\\
$g$ & $0$ & $g$ & $g$ & $c$ & $f$ & $g$ & $1$\\
$1$ & $0$ & $1$ & $1$ & $c$ & $f$ & $1$ & $1$\\
\end{tabular}
\end{table}
\end{example}


\begin{theorem}\label{th3.3}
Let $(L,\leq,0,1)$ be a bounded lattice with $e \in L\setminus \{0,1\}$.\\
$(i)$ Let $T_{e}$ be a $t$-norm on $[0,e]$ and $R$ be a $t$-subconorm on $L$. If $x\parallel y$ for all $x\in I_{e}$ and $y\in [e,1[$, then the function $U_3:L^{2}\rightarrow L$ defined by

$U_{3}(x,y)=\begin{cases}
T_{e}(x, y) &\mbox{if } (x,y)\in [0,e]^{2},\\
x &\mbox{if } (x,y)\in I_{e}\times ]0,e]\cup I_{e}\times ]e,1[,\\
y &\mbox{if } (x,y)\in ]0,e]\times I_{e}\cup ]e,1[\times I_{e},\\
R(x,y) &\mbox{if } (x,y)\in I_{e}^{2}\cup ]e,1]^2,\\
x\vee y &\mbox{if } (x,y)\in D(e)',\\
x\wedge y &\mbox otherwise,
\end{cases}$\\
is a uninorm on $L$  with $e\in L\setminus\{0,1\}$ iff $T_{e}(x,y)>0$ for all $x,y>0$ and $R(x,y)<1$ for all $x,y<1$.\\
$(ii)$ Let $F$ be a $t$-subnorm on $L$ and $S_{e}$ be a $t$-conorm on $[e,1]$. If $x\parallel y$  for all $x\in I_{e}$ and $y\in ]0,e]$, then the function $U_4:L^{2}\rightarrow L$ defined by

$U_{4}(x,y)=\begin{cases}
S_{e}(x, y) &\mbox{if } (x,y)\in [e,1]^{2},\\
x &\mbox{if } (x,y)\in I_{e}\times [e,1[\cup I_{e}\times ]0,e[,\\
y &\mbox{if } (x,y)\in [e,1[\times I_{e}\cup ]0,e[\times I_{e},\\
F(x,y) &\mbox{if } (x,y)\in I_{e}^{2}\cup [0,e[^2,\\
x\wedge y &\mbox{if } (x,y)\in E(e)',\\
x\vee y &\mbox otherwise,
\end{cases}$\\
is a uninorm on $L$  with $e\in L\setminus\{0,1\}$ iff $F(x,y)>0$ for all $x,y>0$ and $S_{e}(x,y)<1$ for all $x,y<1$.

\end{theorem}
\proof Let $(L,\leq,0,1)$ be a bounded lattice with $e \in L\setminus \{0,1\}$. If $x\parallel y$  for all $x\in I_{e}$ and $y\in [e,1[$, then we give the proof of the fact that $U_{3}$ is a uninorm iff $T_{e}(x,y)>0$ for all $x,y>0$ and $R(x,y)<1$ for all $x,y<1$. The same result for $U_{4}$ can be obtained by similar arguments.

Necessity. Let the function $U_{3}$ be a uninorm on $L$ with the identity $e$. Then we will prove that $T_{e}(x,y)>0$ for all $x,y>0$ and $R(x,y)<1$ for all $x,y<1$. Assume that there are some elements $x\in ]0,e[$ and $y\in ]0,e[$ such that $T_{e}(x,y)=0$. If $z\in I_{e}$, then we obtain $U_{3}(x,U_{3}(y,z))=U_{3}(x,z)=z$ and $U_{3}(U_{3}(x,y),z)=U_{3}(T_{e}(x,y),z)=U_{3}(0,z)=0$. Since $T_{e}(x,y)=0$, the associativity property is violated. Assume that there are some elements $x\in ]e,1[$ and $y\in ]e,1[$ such that $R(x,y)=1$. If $z\in I_{e}$, then we obtain $U_{3}(x,U_{3}(y, z))=U_{3}(x,z)=z$ and $U_{3}(U_{3}(x,y), z)=U_{3}(R(x,y),z)=U_{3}(1,z)=1$. Since $R(x,y)=1$, the associativity property is violated. Hence,  $U_{3}$ is not a uninorm on $L$ which is a contradiction. Therefore, the conditions that $T_{e}(x,y)>0$ for all $x,y>0$ and $R(x,y)<1$ for all $x,y<1$ are necessary.

Sufficiency.
To prove the sufficiency, we compare $U_{3}$ with $U_{R}^{e}$ in  Theorem $\ref{th2.2}$  and find that besides the condition that $T_{e}(x,y)>0$ for all $x,y>0$, the values of $U_{3}$ and $U_{R}^{e}$ are different on $I_{e} \times \{0\}\cup \{1\}\times \{0\}$, $\{0\}\times I_{e}\cup \{0\} \times\{1\}$ and $]0,e[\times ]0,e[$. So, combined with  Theorem $\ref{th2.2}$, the proof of the sufficiency  can be reduced to the above areas.

The commutativity of $U_{3}$  and the fact that $e$  is a neutral element are evident. Hence, we prove only the monotonicity and the associativity of $U_{3}$.

I. Monotonicity: We prove that if $x\leq y$, then for all $z\in L$, $U_{3}(x,z)\leq U_{3}(y,z)$.

1.$x=0$ or $y=0$ or $z=0$

\ \ \ \ \ \ \ \ \ \ \ \ $U_{3}(x,z)=0\leq U_{3}(y,z)$

II. Associativity: We need to prove that $U_{3}(x,U_{3}(y,z))=U_{3}(U_{3}(x,y),z)$  for all $x,y,z\in L$.

If one of $x,y,z$ is equal to $0$, then it is clear that $U_{3}(x,U_{3}(y,z))=U_{3}(U_{3}(x,y),z)$.
Hence, $U_{3}$ is a uninorm on $L$ with a neutral element $e$.\\

%
%
%
%
%
%
%
%

It is worth pointing out that if the condition that $x\parallel y$  for all $x\in I_{e}$  and $y\in [e,1[$  in Theorem \ref{th3.3} is not satisfied, then $U_{3}$  may not be a uninorm on $L_{1}$. Similarly, if the condition that $x\parallel y$  for all $x\in I_{e}$  and $y\in ]0,e]$  in Theorem \ref{th3.3} is not satisfied, then $U_{4}$  may not be a uninorm on $L_{1}$.\\

\begin{example}\label{ex3.2}
Consider the bounded lattice $L_{1}=\{0,a,e,c,f,g,1\}$ depicted by the Hasse diagram in Fig.$1$. 

(1) If we take $T_{e}(x,y)=T_{\wedge}(x,y)$  on  $[0,e]$ and $R(x,y)=x\vee y$, then the function $U_{3}$  on $L_{1}$, shown in Table \ref{Tab:03}, is not a uninorm on $L_{1}$  with the neutral element $e$. In fact, the associativity is not satisfied, since $U_{3}(a,U_{3}(c,f))=U_{3}(a,g)=a$ and $U_{3}(U_{3}(a,c),f)=U_{3}(c,f)=g$.

 (2) If we take $S_{e}(x,y)=S_{\vee}(x,y)$  on  $[e,1]$ and $F(x,y)=x\wedge y$, then the function $U_{4}$  on $L_{1}$, shown in Table \ref{Tab:4}, is not a uninorm on $L_{1}$  with the neutral element $e$. In fact, the associativity is not satisfied, since $U_{4}(g,U_{4}(c,f))=U_{4}(g,a)=g$ and $U_{4}(U_{4}(g,c),f)=U_{4}(c,f)=a$.

\begin{table}[htpp]
\centering
\caption{The function $U_{3}$  on $L_{1}$  given in Fig.1.}
\ \ \ \ \\
\label{Tab:03}
\begin{tabular}{c|ccccccc}

 $U_{3}$ & $0$ & $a$ & $e$ & $c$ & $f$ & $g$ & $1$\\
 \hline
$0$ & $0$ & $0$ & $0$ & $0$ & $0$ & $0$ & $0$\\
$a$ & $0$ & $a$ & $a$ & $c$ & $f$ & $a$ & $1$\\
$e$ & $0$ & $a$ & $e$ & $c$ & $f$ & $g$ & $1$\\
$c$ & $0$ & $c$ & $c$ & $c$ & $g$ & $c$ & $1$\\
$f$ & $0$ & $f$ & $f$ & $g$ & $f$ & $f$ & $1$\\
$g$ & $0$ & $a$ & $g$ & $c$ & $f$ & $g$ & $1$\\
$1$ & $0$ & $1$ & $1$ & $1$ & $1$ & $1$ & $1$\\
\end{tabular}
\end{table}
\end{example}


\begin{table}[htpp]
\centering
\caption{The function $U_{4}$  on $L_{1}$  given in Fig.1.}
\ \ \ \ \\
\label{Tab:4}
\begin{tabular}{c|ccccccc}

 $U_{4}$ & $0$ & $a$ & $e$ & $c$ & $f$ & $g$ & $1$\\
 \hline
$0$ & $0$ & $0$ & $0$ & $0$ & $0$ & $0$ & $1$\\
$a$ & $0$ & $a$ & $a$ & $c$ & $f$ & $g$ & $1$\\
$e$ & $0$ & $a$ & $e$ & $c$ & $f$ & $g$ & $1$\\
$c$ & $0$ & $c$ & $c$ & $c$ & $a$ & $c$ & $1$\\
$f$ & $0$ & $f$ & $f$ & $a$ & $f$ & $f$ & $1$\\
$g$ & $0$ & $g$ & $g$ & $c$ & $f$ & $g$ & $1$\\
$1$ & $1$ & $1$ & $1$ & $1$ & $1$ & $1$ & $1$\\
\end{tabular}
\end{table}

\begin{theorem}\label{th3.5}
Let $(L,\leq,0,1)$ be a bounded lattice with $e \in L\setminus \{0,1\}$.\\
$(i)$ Let $T_{e}$ be a $t$-norm on $[0,e]$, and $R$ be a $t$-subconorm on $L$. If $x\parallel y$  for all $x\in I_{e}$  and $y\in[e,1[$,  then the function $U_5:L^{2}\rightarrow L$ defined by

$U_{5}(x,y)=\begin{cases}
T_{e}(x, y) &\mbox{if } (x,y)\in [0,e]^{2},\\
x\wedge y &\mbox{if } (x,y)\in [0,e[\times I_{e}\cup I_{e}\times [0,e[ \cup [0,e[\times ]e,1]\cup ]e,1]\times [0,e[,\\
y &\mbox{if } (x,y)\in \{e\}\times I_{e}\cup \{e\}\times ]e,1],\\
x &\mbox{if } (x,y)\in I_{e}\times \{e\}\cup ]e,1]\times \{e\},\\
R(x,y) &\mbox{if } (x,y)\in {I_{e}}^{2},\\
1 &\mbox{if } (x,y)\in I_{e}\times ]e,1]\cup ]e,1]\times I_{e}\cup ]e,1]^{2},
\end{cases}$\\
 is a uninorm on  $L$ with the neutral element $e$ iff $x> y$  for all $x\in I_{e}$  and $y\in[0,e[$.\\
$(ii)$ Let $S_{e}$ be a $t$-conorm on $[e,1]$, and $F$ be a $t$-subnorm on $L$. If $x\parallel y$  for all $x\in I_{e}$ and $y\in ]0,e]$, then the function $U_{6}:L_{2}\rightarrow L$  defined by

$U_{6}(x,y)=\begin{cases}
S_{e}(x, y) &\mbox{if } (x,y)\in [e,1]^{2},\\
x\vee y &\mbox{if } (x,y)\in ]e,1]\times I_{e}\cup I_{e}\times ]e,1] \cup [0,e[\times ]e,1]\cup ]e,1]\times [0,e[,\\
y &\mbox{if } (x,y)\in \{e\}\times I_{e}\cup \{e\}\times [0,e[,\\
x &\mbox{if } (x,y)\in I_{e}\times \{e\}\cup [0,e[\times \{e\},\\
F(x,y) &\mbox{if } (x,y)\in {I_{e}}^{2},\\
0 &\mbox{if } (x,y)\in I_{e}\times [0,e[\cup [0,e[\times I_{e}\cup [0,e[^{2},
\end{cases}$\\
is a uninorm on $L$  with the neutral element $e$ iff $x<y$   for all $x\in I_{e}$ and $y\in ]e,1]$.
\end{theorem}
\proof
Let $(L,\leq,0,1)$ be a bounded lattice with $e \in L\setminus \{0,1\}$. If $x\parallel y$ for all $x\in I_{e}$ and $y\in [e,1[$, then we give the proof of the fact that $U_{5}$ is a uninorm iff $x> y$  for all $x\in I_{e}$  and $y\in[0,e[$. The same result for $U_{6}$ can be obtained using similar arguments.

Necessity. Let the function $U_{5}$ be a uninorm on $L$ with the identity $e$. Then we will prove that $x>y$ for all $x\in I_{e}$ and $y\in [0,e[$.  Assume that there are some elements $x\in I_{e}$ and $y\in [0,e[$ such that $x\parallel y$. If $z\in ]e,1]$, then  $U_{5}(y, U_{5}(x, z))=U_{5}(y, 1)=y$ and $U_{5}(U_{5}(y, x), z)=U_{5}(x\wedge y, z)=x\wedge y\wedge z=x\wedge y$. Since $x\parallel y$, the associativity property is violated. Then $U_{5}$ is not a uninorm on $L$ which is a contradiction. Hence, $x>y$ for all $x\in I_{e}$ and $y\in [0,e[$.

Sufficiency. Let $x>y$ for all $x\in I_{e}$ and $y\in [0,e[$, with $x\parallel y$  for all $x\in I_{e}$  and $y\in[e,1[$. The commutativity of $U_{5}$  and the fact that $e$  is a neutral element are evident. Hence, we prove only the monotonicity and the associativity of $U_{5}$.

I. Monotonicity: We prove that if $x\leq y$, then for all $z\in L$, $U_{5}(x,z)\leq U_{5}(y,z)$. The proof is split into all possible cases.
If one of $x$, $y$, $z$ is equal to $e$, then it is clear that $U_{5}(x,z)=z=U_{5}(y,z)$

1. Let $z=e$.

\ \ \ \ \ \ \ \ \ \ \ \ $U_5(x,z)=x\leq y=U_5(y,z)$

2. Let $x=e$.

\ \ \ 2.1. $y=e$,

\ \ \ \ \ \ \ \ \ \ \ \ $U_5(x,z)=z=U_5(y,z)$

\ \ \ 2.2. $y\in ]e,1]$,

\ \ \ \ \ \ 2.2.1. $z\in [0,e[$,

\ \ \ \ \ \ \ \ \ \ \ \ $U_5(x,z)=z=U_5(y,z)$

\ \ \ \ \ \ 2.2.2. $z\in ]e,1]$,

\ \ \ \ \ \ \ \ \ \ \ \ $U_5(x,z)=z\leq 1=U_5(y,z)$

\ \ \ \ \ \ 2.2.3. $z\in I_e$,

\ \ \ \ \ \ \ \ \ \ \ \ $U_5(x,z)=z\leq 1=U_5(y,z)$

3. Let $y=e$.

\ \ \ 3.1. $x\in [0,e[$,

\ \ \ \ \ \ 3.1.1. $z\in [0,e[$,

\ \ \ \ \ \ \ \ \ \ \ \ $U_5(x,z)=T_e(x,z)\leq z=U_5(y,z)$

\ \ \ \ \ \ 3.1.2. $z\in ]e,1]$,

\ \ \ \ \ \ \ \ \ \ \ \ $U_5(x,z)=x\leq z=U_5(y,z)$

\ \ \ \ \ \ 3.1.3. $z\in I_e$,

\ \ \ \ \ \ \ \ \ \ \ \ $U_5(x,z)=x\wedge z=x\leq z=U_5(y,z)$

There is no $e$ in any of the following cases.

1. Let $x\in [0,e[$.

\ \ \ 1.1. $y\in [0,e[$,

\ \ \ \ \ \ 1.1.1. $z\in [0,e[$,

\ \ \ \ \ \ \ \ \ \ \ \ $U_5(x,z)=T_e(x,z)\leq T_e(y,z)=U_5(y,z)$

\ \ \ \ \ \ 1.1.2. $z\in ]e,1]$,

\ \ \ \ \ \ \ \ \ \ \ \ $U_5(x,z)=x\wedge z=x\leq y= y\wedge z=U_5(y,z)$

\ \ \ \ \ \ 1.1.3. $z\in I_e$,

\ \ \ \ \ \ \ \ \ \ \ \ $U_5(x,z)=x\wedge z=x\leq y= y\wedge z=U_5(y,z)$

\ \ \ 1.2. $y\in ]e,1]$,

\ \ \ \ \ \ 1.2.1. $z\in [0,e[$,

\ \ \ \ \ \ \ \ \ \ \ \ $U_5(x,z)=T_e(x,z)\leq z= y\wedge z=U_5(y,z)$

\ \ \ \ \ \ 1.2.2. $z\in ]e,1]$,

\ \ \ \ \ \ \ \ \ \ \ \ $U_5(x,z)=x\wedge z=x\leq 1=U_5(y,z)$

\ \ \ \ \ \ 1.2.3. $z\in I_e$,

\ \ \ \ \ \ \ \ \ \ \ \ $U_5(x,z)=x\wedge z=x\leq 1=U_5(y,z)$

\ \ \ 1.3. $y\in I_e$,

\ \ \ \ \ \ 1.3.1. $z\in [0,e[$,

\ \ \ \ \ \ \ \ \ \ \ \ $U_5(x,z)=T_e(x,z)\leq z= y\wedge z=U_5(y,z)$

\ \ \ \ \ \ 1.3.2. $z\in ]e,1]$,

\ \ \ \ \ \ \ \ \ \ \ \ $U_5(x,z)=x\wedge z=x\leq 1=U_5(y,z)$

\ \ \ \ \ \ 1.3.3. $z\in I_e$,

\ \ \ \ \ \ \ \ \ \ \ \ $U_5(x,z)=x\wedge z=x\leq R(y,z)=U_5(y,z)$

2. Let $x\in ]e,1]$.

\ \ \ 2.1. $y\in ]e,1]$,

\ \ \ \ \ \ 2.1.1. $z\in [0,e[$,

\ \ \ \ \ \ \ \ \ \ \ \ $U_5(x,z)=x\wedge z=z=y\wedge z=U_5(y,z)$

\ \ \ \ \ \ 2.1.2. $z\in ]e,1]$,

\ \ \ \ \ \ \ \ \ \ \ \ $U_5(x,z)=1=U_5(y,z)$

\ \ \ \ \ \ 2.1.3. $z\in I_e$,

\ \ \ \ \ \ \ \ \ \ \ \ $U_5(x,z)=1=U_5(y,z)$

3. Let $x\in I_e$.

\ \ \ 3.1. $y\in I_e$,

\ \ \ \ \ \ 3.1.1. $z\in [0,e[$,

\ \ \ \ \ \ \ \ \ \ \ \ $U_5(x,z)=x\wedge z=z=y\wedge z=U_5(y,z)$

\ \ \ \ \ \ 3.1.2. $z\in ]e,1]$,

\ \ \ \ \ \ \ \ \ \ \ \ $U_5(x,z)=1=U_5(y,z)$

\ \ \ \ \ \ 3.1.3. $z\in I_e$,

\ \ \ \ \ \ \ \ \ \ \ \ $U_5(x,z)=R(x,z)\leq R(y,z)=U_5(y,z)$

II. Associativity: We need to prove that $U_5(x,U_5(y,z))=U_5(U_5(x,y),z)$  for all $x,y,z\in L$. The proof is split into all possible cases.
If one of $x,y,z$ is equal to $e$, then it is clear that $U_5(x,U_5(y,z))=U_5(U_5(x,y),z)$. So there is no $e$ in any of the following cases. By Proposition  \ref{pro2.2}, we only need to consider the following cases:

1. Let $x\in [0,e[$.

\ \ \ 1.1. $y\in [0,e[$,

\ \ \ \ \ \ 1.1.1. $z\in [0,e[$,

\ \ \ \ \ \ \ \ \ \ \ \ $U_5(x,U_5(y,z))=T_e(x,T_e(y,z))=T_e(T_e(x,y),z)=U_5(U_5(x,y),z)$

\ \ \ \ \ \ 1.1.2. $z\in ]e,1]$,

\ \ \ \ \ \ \ \ \ \ \ \ $U_5(x,U_5(y,z))=U_5(x,y\wedge z)=T_e(x,y)=U_5(T_e(x,y),z)=$

\ \ \ \ \ \ \ \ \ \ \ \ $U_5(U_5(x,y),z)$

\ \ \ \ \ \ 1.1.3. $z\in I_e$,

\ \ \ \ \ \ \ \ \ \ \ \ $U_5(x,U_5(y,z))=U_5(x,y\wedge z)=T_e(x,y)=U_5(T_e(x,y),z)=$

\ \ \ \ \ \ \ \ \ \ \ \ $U_5(U_5(x,y),z)$

\ \ \ 1.2. $y\in ]e,1]$,

\ \ \ \ \ \ 1.2.1. $z\in]e,1]$,

\ \ \ \ \ \ \ \ \ \ \ \ $U_{5}(x, U_{5}(y, z))=U_{5}(x, 1)=x=U_{5}(x, z)=U_{5}(U_{5}(x, y), z)$

\ \ \ \ \ \ 1.2.2. $z\in I_{e}$,

\ \ \ \ \ \ \ \ \ \ \ \ $U_{5}(x, U_{5}(y, z))=U_{5}(x, 1)=x, U_{5}(U_{5}(x, y), z)=U_{5}(x, z)=x$

\ \ \ \ \ \ \ \ \ \ \ \ and $U_{5}(y, U_{5}(x, z))=U_{5}(y, x)=x$

\ \ \ 1.3. $y\in I_{e}$,

\ \ \ \ \ \ 1.3.1. $z \in I_{e}$,

\ \ \ \ \ \ \ \ \ \ \ \ $U_{5}(x, U_{5}(y, z))=U_{5}(x, R(y, z))=x=U_{5}(x, z)=U_{5}(U_{5}(x, y), z)$

2. Let $x\in] e, 1]$.

\ \ \ 2.1. $y\in]e, 1]$,

\ \ \ \ \ \ 2.1.1. $z\in] e, 1]$,

\ \ \ \ \ \ \ \ \ \ \ \ $U_{5}(x,U_{5}(y,z))=U_{5}(x, 1)=1=U_{5}(1, z)=U_{5}(U_{5}(x, y), z)$

\ \ \ \ \ \ 2.1.2. $z\in I_{e}$,

\ \ \ \ \ \ \ \ \ \ \ \ $U_{5}(x, U_{5}(y, z))=U_{5}(x, 1)=1=U_{5}(1, z)=U_{5}(U_{5}(x, y), z)$

\ \ \ 2.2. $y\in I_{e}$,

\ \ \ \ \ \ 2.2.1. $z \in I_{e}$,

\ \ \ \ \ \ \ \ \ \ \ \ $U_{5}(x, U_{5}(y, z))=U_{5}(x, R(y, z))=1=U_{5}(1, z)=U_{5}(U_{5}(x, y), z)$

3. Let $x \in I_{e}$.

\ \ \ 3.1. $y \in I_{e}$,

\ \ \ \ \ \ 3.1.1. $z \in I_{e}$,

\ \ \ \ \ \ \ \ \ 3.1.1.1. $R(x, y)=1, R(y, z)=1$,

\ \ \ \ \ \ \ \ \ \ \ \ $U_{5}(x, U_{5}(y, z))=U_{5}(x, R(y, z))=1=U_{5}(R(x, y), z)=U_{5}(U_{5}(x, y), z)$

\ \ \ \ \ \ \ \ \ 3.1.1.2. $R(x, y)=1, R(y, z) \in I_{e}$,
\begin{eqnarray*} U_{5}(x, U_{5}(y, z))&=&U_{5}(x, R(y, z))=R(x, R(y, z))\\
&=&R(R(x, y), z)=
R(1, z)=1\\&=&U_{5}(1, z)=U_{5}(R(x, y), z)\\
&=&U_{5}(U_{5}(x, y),z) (\mathrm{by\ Remark\ \ref{Re2.8}(3)})
\end{eqnarray*}

\ \ \ \ \ \ \ \ \ 3.1.1.3. $R(x, y) \in I_{e}, R(y, z)=1$,
\begin{eqnarray*}
U_{5}(x, U_{5}(y, z)) &=& U_{5}(x, R(y, z))=U_{5}(x, 1)\\
&=&1=R(x, 1)=R(x, R(y, z))\\
&=&R(R(x, y), z)=U_{5}(R(x, y), z)\\
&=&U_{5}(U_{5}(x, y), z) (\mathrm{by\ Remark\ \ref{Re2.8}(3)})
\end{eqnarray*}

\ \ \ \ \ \ \ \ \ 3.1.1.4. $R(x, y) \in I_{e}, R(y, z) \in I_{e}$,

\ \ \ \ \ \ \ \ \ \ \ \ $U_{5}(x, U_{5}(y, z))=R(x, R(y, z))=R(R(x, y), z)=U_{5}(U_{5}(x, y), z)$

Hence, we obtain that $U_{5}(x, U_{5}(y, z))=U_{5}(U_{5}(x, y), z)$ for all $x,y,z\in L$ by Proposition \ref{pro2.2}. Therefore, $U_{5}$ is a uninorm on $L$ with the neutral element $e$.\\

It is worth pointing out that if the condition that $x\parallel y$  for all $x\in I_{e}$  and $y\in [e,1[$  in Theorem \ref{th3.5} is not satisfied, then $U_{5}$  may not be a uninorm on $L_{2}$  with the neutral element $e$. Similarly, if the condition that $x\parallel y$  for all $x\in I_{e}$  and $y\in ]0,e]$  in Theorem \ref{th3.5} is not satisfied, then $U_{6}$  may not be a uninorm on $L_{2}$  with the neutral element $e$. The following example will show the above facts.

\begin{example}\label{ex3.6}
Consider the bounded lattice  $L_{2}=\{0,a,e,c,d,f,g,1\}$  depicted by the Hasse diagram in Fig.$2$, where $a\in [0,e[$, $g\in ]e,1]$, $c\in I_{e}$, $d\in I_{e}$, $f\in I_{e}$, $g>c$, $g>d$, $g>f$, $a<c$, $a<d$ and $a<f$.
If we take $T_{e}(x,y)=T_{\wedge}(x,y)$  on  $[0,e]$ and $R(x,y)=x\vee y$, then  the function $U_{5}$  on $L_{2}$, shown in Table \ref{Tab:05}, is not a uninorm on $L_{2}$  with the neutral element $e$. In fact, the associativity is not satisfied, since $U_{5}(f,U_{5}(c,d))=U_{5}(f,d)=g$ and $U_{5}(U_{5}(f,c),d)=U_{5}(g,d)=1$.
 Similarly, if we take $S_e(x,y)=S_{\vee}(x,y)$ on $[e,1]$ and  $F(x,y)=x\wedge y$, then  the function $U_{6}$  on $L_{2}$, shown in Table \ref{Tab:06}, is not a uninorm on $L_{2}$  with the neutral element $e$. In fact, the associativity is not satisfied, since $U_{6}(f,U_{6}(c,d))=U_{6}(f,c)=a$ and $U_{6}(U_{6}(f,c),d)=U_{6}(a,d)=0$.
\par\noindent\vskip50pt
\begin{minipage}{11pc}
\setlength{\unitlength}{0.75pt}
\begin{picture}(600,180)
\put(240,40){\circle{5}}\put(238,29){\makebox(0,0)[l]{\footnotesize$0$}}
\put(240,80){\circle{5}}\put(245,80){\makebox(0,0)[l]{\footnotesize$a$}}
\put(200,120){\circle{5}}\put(189,120){\makebox(0,0)[l]{\footnotesize$c$}}
\put(220,140){\circle{5}}\put(210,145){\makebox(0,0)[l]{\footnotesize$d$}}
\put(240,120){\circle{5}}\put(245,120){\makebox(0,0)[l]{\footnotesize$e$}}
\put(280,120){\circle{5}}\put(285,120){\makebox(0,0)[l]{\footnotesize$f$}}
\put(240,160){\circle{5}}\put(245,164){\makebox(0,0)[l]{\footnotesize$g$}}
\put(240,200){\circle{5}}\put(237,213){\makebox(0,0)[l]{\footnotesize$1$}}
\put(240,43){\line(0,1){34}}
\put(240,83){\line(0,1){34}}
\put(240,123){\line(0,1){34}}
\put(240,163){\line(0,1){34}}
\put(238,82){\line(-1,1){36}}
\put(242,82){\line(1,1){36}}
\put(202,122){\line(1,1){16}}
\put(222,142){\line(1,1){16}}
\put(278,122){\line(-1,1){36}}
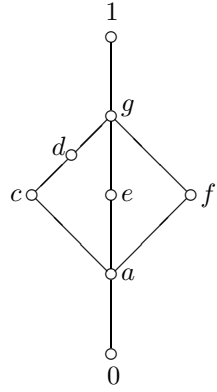
\put(170,-10){\emph{Fig.2. The lattice $L_{2}$}}
\end{picture}
\end{minipage}

\begin{table}[htpp]
\centering
\caption{The function $U_{5}$  on $L_{2}$  given in Fig.2.}
\ \ \ \ \\
\label{Tab:05}
\begin{tabular}{c|cccccccc}

 $U_{5}$ & $0$ & $a$ & $e$ & $c$ & $d$ & $f$ & $g$ & $1$\\
 \hline
$0$ & $0$ & $0$ & $0$ & $0$ & $0$ & $0$ & $0$ & $0$\\
$a$ & $0$ & $a$ & $a$ & $a$ & $a$ & $a$ & $a$ & $a$\\
$e$ & $0$ & $a$ & $e$ & $c$ & $d$ & $f$ & $g$ & $1$\\
$c$ & $0$ & $a$ & $c$ & $c$ & $d$ & $g$ & $1$ & $1$\\
$d$ & $0$ & $a$ & $d$ & $d$ & $d$ & $g$ & $1$ & $1$\\
$f$ & $0$ & $a$ & $f$ & $g$ & $g$ & $f$ & $1$ & $1$\\
$g$ & $0$ & $a$ & $g$ & $1$ & $1$ & $1$ & $1$ & $1$\\
$1$ & $0$ & $a$ & $1$ & $1$ & $1$ & $1$ & $1$ & $1$\\
\end{tabular}
\end{table}

\par\noindent\vskip50pt

\begin{table}[htpp]
\centering
\caption{The function $U_{6}$  on $L_{2}$  given in Fig.2.}
\ \ \ \ \\
\label{Tab:06}
\begin{tabular}{c|cccccccc}

 $U_{6}$ & $0$ & $a$ & $e$ & $c$ & $d$ & $f$ & $g$ & $1$\\
 \hline
$0$ & $0$ & $0$ & $0$ & $0$ & $0$ & $0$ & $g$ & $1$\\
$a$ & $0$ & $0$ & $a$ & $0$ & $0$ & $0$ & $g$ & $1$\\
$e$ & $0$ & $a$ & $e$ & $c$ & $d$ & $f$ & $g$ & $1$\\
$c$ & $0$ & $0$ & $c$ & $c$ & $c$ & $a$ & $g$ & $1$\\
$d$ & $0$ & $0$ & $d$ & $c$ & $d$ & $a$ & $g$ & $1$\\
$f$ & $0$ & $0$ & $f$ & $a$ & $a$ & $f$ & $g$ & $1$\\
$g$ & $g$ & $g$ & $g$ & $g$ & $g$ & $g$ & $g$ & $1$\\
$1$ & $1$ & $1$ & $1$ & $1$ & $1$ & $1$ & $1$ & $1$\\
\end{tabular}
\end{table}
\end{example}

\begin{corollary}\label{co3.1}
Let $(L,\leq,0,1)$ be a bounded lattice with $e\in L\setminus \{0,1\}$,  $x\parallel y$  for all $x\in I_{e}$ and $y\in [e,1[$, $x>y$  for all  $x\in I_{e}$ and $y\in [0,e[$. Suppose that $U_{5}$  is the uninorm on $L$ defined as in Theorem \ref{th3.5},   $U_{R}$ is the uninorm in Theorem \ref{th2.3} and $U_{\wedge}$  is the uninorm  in Theorem \ref{th2.4}. Then $U_{5}\geq U_{R}$ and  $U_{5}\geq U_{\wedge}$.
\end{corollary}


\begin{corollary}\label{co3.2}
Let $(L,\leq,0,1)$  be a bounded lattice with $e\in L\setminus \{0,1\}$, $x\parallel y$  for all $x\in I_{e}$ and $y\in ]0,e]$,  $x<y$  for all $x\in I_{e}$ and $y\in ]e,1]$. Suppose that $U_{6}$  is the uninorm on $L$ defined as in Theorem \ref{th3.5}, $U_{F}$  is the uninorm in Theorem \ref{th2.3} and $U_{\vee}$  is the uninorm   in Theorem \ref{th2.4}. Then $U_{6}\leq U_{F}$ and  $U_{6}\leq U_{\vee}$.
\end{corollary}

\section{Modifications of some  uninorms}

In this section, we mainly modify some  uninorms, in which some constraints  are actually necessary and sufficient.

\begin{theorem}\label{th3.2}
Let $(L,\leq,0,1)$ be a bounded lattice with $e \in L\setminus \{0,1\}$.\\
(i) Let $T_{e}$ be a $t$-norm on $[0,e]$ and $S_{e}$ be a $t$-conorm on $[e,1]$. If $x\parallel y$ for all $x\in I_{e}$ and $y\in [e,1[$, then the function $U_{1}^{e}:L^{2}\rightarrow L$ defined by

$U_{1}^{e}(x,y)=\begin{cases}
T_{e}(x, y) &\mbox{if } (x,y)\in [0,e]^{2},\\
S_{e}(x, y) &\mbox{if } (x,y)\in [e,1]^{2},\\
x &\mbox{if } (x,y)\in I_{e}\times [e,1[\cup I_{e}\times ]0,e[,\\
y &\mbox{if } (x,y)\in [e,1[\times I_{e}\cup ]0,e[\times I_{e},\\
x\vee y &\mbox{if } (x,y)\in I_{e}^2\cup I_{e}\times \{1\}\cup \{1\}\times I_{e}\cup ]0,e[\times \{1\}\cup \{1\}\times ]0,e[,\\
x\wedge y &\mbox{  } otherwise,\\
\end{cases}$\\
is a uninorm on $L$  with $e\in L\setminus\{0,1\}$ iff $T_{e}(x,y)>0$ for all $x,y>0$ and $S_{e}(x,y)<1$ for all $x,y<1$.\\
(ii) Let $T_{e}$ be a $t$-norm on $[0,e]$ and $S_{e}$ be a $t$-conorm on $[e,1]$. If $x\parallel y$  for all $x\in I_{e}$ and $y\in ]0,e]$, then the function $U_{2}^{e}:L^{2}\rightarrow L$ defined by

$U_{2}^{e}(x,y)=\begin{cases}
T_{e}(x, y) &\mbox{if } (x,y)\in [0,e]^{2},\\
S_{e}(x, y) &\mbox{if } (x,y)\in [e,1]^{2},\\
x &\mbox{if } (x,y)\in I_{e}\times ]e,1[\cup I_{e}\times ]0,e],\\
y &\mbox{if } (x,y)\in ]e,1[\times I_{e}\cup ]0,e]\times I_{e},\\
x\wedge y &\mbox{if } (x,y)\in I_{e}^2\cup I_{e}\times \{0\}\cup \{0\}\times I_{e}\cup ]e,1[\times \{0\}\cup \{0\}\times ]e,1[,\\
x\vee y &\mbox{  } otherwise,\\
\end{cases}$\\
is a uninorm on $L$  with $e\in L\setminus\{0,1\}$ iff $T_{e}(x,y)>0$ for all $x,y>0$ and $S_{e}(x,y)<1$ for all $x,y<1$.

\end{theorem}

\proof
If $x\parallel y$  for all $x\in I_{e}$ and $y\in [e,1[$, then we give the proof of the fact that $U_{1}^{e}$ is a uninorm iff $T_{e}(x,y)>0$ for all $x,y> 0$ and $S_{e}(x,y)<1$ for all $x,y<1$. The same result for $U_{2}^{e}$ can be obtained by similar arguments.

Necessity. Let $U_{1}^{e}$ be a uninorm on $L$ with the identity $e$. Next we will  prove that $T_{e}(x,y)>0$ for all $x,y>0$ and $S_{e}(x,y)<1$ for all $x,y<1$.
Assume that there exist $x\in ]0,e[$ and $y\in ]0,e[$ such that $T_{e}(x,y)=0$. Take $z\in I_{e}$.   Then  $U_{1}^{e}(x,U_{1}^{e}(y,z))=U_{1}^{e}(x,z)=z$ and, since $T_{e}(x,y)=0$, $U_{1}^{e}(U_{1}^{e}(x,y),z)=U_{1}^{e}(T_{e}(x,y),z)=U_{1}^{e}(0,z)=0$. Hence, the associativity property  is violated. Assume that there exist $x\in ]e,1[$ and $y\in ]e,1[$ such that $S_{e}(x,y)=1$. Take  $z\in I_{e}$. Then  $U_{1}^{e}(x,U_{1}^{e}(y, z))=U_{1}^{e}(x,z)=z$ and,  since $S_{e}(x,y)=1$, $U_{1}^{e}(U_{1}^{e}(x,y), z)=U_{1}^{e}(S_{e}(x,y),z)=U_{1}^{e}(1,z)=1$.  Hence, the associativity property is violated. So,  $U_{1}^{e}$ is not a uninorm on $L$ which is a contradiction. Therefore, the conditions that $T_{e}(x,y)>0$ for all $x,y>0$ and $S_{e}(x,y)<1$ for all $x,y<1$ are necessary.

Sufficiency. This proof is the same as Theorem \ref{th2.1}.\\

\begin{remark}\label{411}
In Theorem \ref{th2.1},  the   conditions that   $T_{e}$  is a $t$-norm on $[0,e]$ such that $T_{e}(x,y)>0$ for all $x,y>0$ and $S_{e}$  is a $t$-conorm on $[e,1]$ such that $S_{e}(x,y)<1$ for all $x,y<1$ are sufficient for $U_{1}^{e}$ and $U_{2}^{e}$.   In fact, these conditions are also necessary for $U_{1}^{e}$ and $U_{2}^{e}$ in Theorem \ref{th3.2}.
\end{remark}


\begin{remark}\label{412}
Let $U_{1}$ and $U_{2}$  be the uninorms on $L$ defined as in Theorem \ref{th3.1}, and let $U_{1}^{e}$  and $U_{2}^{e}$ be the uninorms on $L$ defined as in Theorem \ref{th3.2}.

(1) It is easy to see that  the values of  $U_{1}$ and $U_{1}^{e}$ are  different on  $I_{e} \times \{0\}\cup \{1\}\times \{0\}$ and $\{0\}\times I_{e}\cup \{0\} \times\{1\}$.
In fact, $U_{1}=x$ and $U_{1}^{e}=0$ when $(x,y)\in I_{e} \times \{0\}\cup \{1\}\times \{0\}$; $U_{1}=y$, $U_{1}^{e}=0$ when $(x,y)\in \{0\}\times I_{e}\cup \{0\} \times\{1\}$. The values of $U_{1}$ could be $0$ on $]0,e[\times ]0,e[$;  however, the values of $U_{1}^{e}$ are not equal to $0$ on $]0,e[\times ]0,e[$.


(2) The additional constraints of  $U_{1}$ and $U_{1}^{e}$ are  different, that is, the condition that $T_{e}(x,y)>0$ for all $x,y>0$  is not needed for $U_{1}$.

According to (1) and (2),  by changing the  values of  $U_{1}^{e}$, we can  obtain  a new method  to construct the  nuinorm  $U_{1}$  and  then  additional constraints are  different for $U_{1}^{e}$ and  $U_{1}$.

Similarly,  the relationship between  $U_{2}$ and $U_{2}^{e}$   is  same as that of  $U_{1}$ and $U_{1}^{e}$.

\end{remark}


\begin{theorem}\label{th3.4}
Let $(L,\leq,0,1)$ be a bounded lattice with $e \in L\setminus \{0,1\}$.\\
(i) Let $T_{e}$ be a $t$-norm on $[0,e]$ and $R$ be a $t$-subconorm on $L$. If $x\parallel y$ for all $x\in I_{e}$ and $y\in [e,1[$, then the function $U_{R}^{e}:L^{2}\rightarrow L$ defined by

$U_{R}^{e}(x,y)=\begin{cases}
T_{e}(x, y) &\mbox{if } (x,y)\in [0,e]^{2},\\
x &\mbox{if } (x,y)\in I_{e}\times [0,e]\cup I_{e}\times [e,1[,\\
y &\mbox{if } (x,y)\in [0,e]\times I_{e}\cup [e,1[\times I_{e},\\
R(x, y) &\mbox{if } (x,y)\in I_{e}^2\cup ]e,1]^{2},\\
x\vee y &\mbox{if } (x,y)\in D(e),\\
x\wedge y &\mbox{  } otherwise,\\
\end{cases}$\\
is a uninorm on $L$  with $e\in L\setminus\{0,1\}$ iff $R(x,y)<1$ for all $x,y<1$.\\
(ii) Let $S_{e}$ be a $t$-conorm on $[e,1]$ and $F$ be a $t$-subnorm on $L$. If $x\parallel y$  for all $x\in I_{e}$ and $y\in ]0,e]$, then the function $U_{F}^{e}:L^{2}\rightarrow L$ defined by

$U_{F}^{e}(x,y)=\begin{cases}
S_{e}(x, y) &\mbox{if } (x,y)\in [e,1]^{2},\\
x &\mbox{if } (x,y)\in I_{e}\times ]0,e]\cup I_{e}\times [e,1],\\
y &\mbox{if } (x,y)\in ]0,e]\times I_{e}\cup [e,1]\times I_{e},\\
F(x, y) &\mbox{if } (x,y)\in I_{e}^2\cup [0,e[^{2},\\
x\wedge y &\mbox{if } (x,y)\in E(e),\\
x\vee y &\mbox{  } otherwise,\\
\end{cases}$\\
is a uninorm on $L$  with $e\in L\setminus\{0,1\}$ iff $F(x,y)>0$ for all $x,y>0$.

\end{theorem}

\proof
If $x\parallel y$ for all $x\in I_{e}$ and $y\in [e,1[$, then we give the proof of the fact that $U_{R}^{e}$ is a uninorm iff $R(x,y)<1$ for all $x,y<1$. The same result for $U_{F}^{e}$ can be obtained using similar arguments.

Necessity. Let   $U_{R}^{e}$ be a uninorm on $L$ with the identity $e$. Then we will  prove that $R(x,y)<1$ for all $x,y<1$. Assume that there exist $x\in ]e,1[$ and $y\in ]e,1[$ such that $R(x,y)=1$. Take $z\in I_{e}$. Then   $U_{R}^{e}(x,U_{R}^{e}(y, z))=U_{R}^{e}(x,z)=z$ and, since $R(x,y)=1$,  $U_{R}^{e}(U_{R}^{e}(x,y), z)=U_{R}^{e}(R(x,y),z)=U_{R}^{e}(1,z)=1$. Hence,  the associativity   property of   $U_{R}^{e}$ is violated. So,   $U_{R}^{e}$ is not a uninorm on $L$ which is a contradiction. Therefore, the condition that $R(x,y)<1$ for all $x,y<1$ is necessary.

Sufficiency. This proof is the same as Theorem \ref{th2.2}.\\

\begin{remark}\label{re3.3}
In Theorem \ref{th2.2}, the conditions that $R$ is a $t$-subconorm on $L$ such that $R(x,y)<1$ for all $x,y<1$ and $F$ is a $t$-subnorm on $[0,e]$ such that $F(x,y)>0$ for all $x,y>0$ are  sufficient for $U_{R}^{e}$ and $U_{F}^{e}$, respectively. In fact, the conditions are also necessary for $U_{F}^{e}$ and $U_{R}^{e}$ in Theorem \ref{th3.4}, respectively.
\end{remark}

\begin{remark}\label{re3.4}
Let  $U_{3}$ and $U_{4}$ be the uninorms on $L$ defined as in Theorem \ref{th3.3}, and let  $U_{R}^{e}$ and $U_{F}^{e}$ be the uninorms on $L$ defined as in Theorem \ref{th3.4}.

(1) It is easy to see that the values of $U_{3}$ and $U_{R}^{e}$ are different on $I_{e} \times \{0\}\cup \{1\}\times \{0\}$ and $\{0\}\times I_{e}\cup \{0\} \times\{1\}$.   In fact, $U_{3}=0$ and $U_{R}^{e}=x$ when $(x,y)\in I_{e} \times \{0\}\cup \{1\}\times \{0\}$;  $U_{3}=0$ and $U_{R}^{e}=y$ when $(x,y)\in \{0\}\times I_{e}\cup \{0\} \times\{1\}$. The values of $U_{3}$ are not equal to $0$ on $]0,e[\times ]0,e[$; however,  the values of $U_{R}^{e}$ could be $0$ on $]0,e[\times ]0,e[$.


(2) The additional constraints of  $U_{3}$ and $U_{R}^{e}$ are  different, that is, the condition that $T_{e}(x,y)>0$ for all $x,y>0$  is  needed for $U_{3}$.

According to (1) and (2),  by changing the  values of   $U_{R}^{e}$, we can  obtain  a new method  to construct the  nuinorm $U_{3}$  and  then  additional constraints are  different for $U_{R}^{e}$ and  $U_{3}$ .

Similarly,  the relationship between  $U_{4}$ and $U_{F}^{e}$   is  same as that of  $U_{3}$ and $U_{R}^{e}$.

%

\end{remark}

\section{ New constructions of uninorms by interation}
  In this section,  we present new constructions of uninorms  on bounded lattices   by interation based on  $t$-conorms  and $t$-conorms.

\begin{theorem}\label{th3.0}
Let $(L,\leq, 0, 1)$ be a bounded lattice.\\
(i) Let $\{a_{0},a_{1},a_{2},\ldots,a_{n}\}$  be a finite chain in $L$ with $0_{L}=a_{0}<a_{1}<a_{2}<\ldots<a_{n}=1_{L}$ on $(L,\leq,0,1)$,
$T_{e}$ be a $t$-norm on $[0,a_{1}]$ and $S_{i}$ be a $t$-conorm on $[a_{i-1},a_{i}]$ of $L$  for $i\in \{2,3,\ldots,n\}$.  Then for $i\in \{2,3,\ldots,n\}$,   the function $U_{i}^{S}:[0,a_{i}]^{2}\rightarrow [0,a_{i}]$ defined as follows

$U_{i}^{S}(x,y)=\begin{cases}
U_{i-1}^{S}(x, y) &\mbox{if } (x,y)\in [0,a_{i-1}]^{2},\\
S_{i}(x, y) &\mbox{if } (x,y)\in ]a_{i-1},a_{i}]^{2},\\
x\vee y &\mbox{if } (x,y)\in [0,e]\times (I_{a_{i-1}}^{e}\cup]a_{i-1},a_{i}])\\
&\mbox{ } \ \ \ \ \ \ \ \ \ \ \ \cup (I_{a_{i-1}}^{e}\cup]a_{i-1},a_{i}])\times[0,e],\\
x &\mbox{if } (x,y)\in I_{e,a_{i-1}}\times [0,e],\\
y &\mbox{if } (x,y)\in [0,e]\times I_{e,a_{i-1}},\\
S_{i}(x\vee a _{i-1}, y\vee a _{i-1}) &\mbox{  }  otherwise,\\
\end{cases}$\\
is a uninorm on $[0,a_{i}]$, where $U_{1}^{S}=T_{e}$.\\
(ii) Let $\{b_{0},b_{1},b_{2},\ldots,b_{n}\}$  be a finite chain in $L$ with $0_{L}=b_{n}\leq b_{n-1}\leq b_{n-2}\leq \ldots \leq b_{0}=1$ on $(L,\leq,0,1)$,
$S_{e}$ be a t-conorm on $[b_{1},1]$ and $T_{i}$ be a t-norm on $[b_{i},b_{i-1}]$ for $i\in \{2,3,\ldots,n\}$. Then  for $i\in \{2,3,\ldots,n\}$,  the function $U_{i}^{T}:[b_{i},1]\rightarrow [b_{i},1]$ defined as follows

$U_{i}^{T}(x,y)=\begin{cases}
U_{i-1}^{T} &\mbox{if } (x,y)\in [b_{i-1},1]^{2},\\
T_{i}(x, y) &\mbox{if } (x,y)\in [b_{i},b_{i-1}[^{2},\\
x\wedge y &\mbox{if } (x,y)\in [e,1]\times (I_{b_{i-1}}^{e}\cup[b_{i},b_{i-1}[)\\
&\mbox{ } \ \ \ \ \ \ \ \ \ \ \ \cup (I_{b_{i-1}}^{e}\cup[a_{i},b_{i-1}[)\times[e,1],\\
x &\mbox{if } (x,y)\in I_{e,b_{i-1}}\times [e,1],\\
y &\mbox{if } (x,y)\in [e,1]\times I_{e,b_{i-1}},\\
T_{i}(x\wedge b_{i-1}, y\wedge b_{i-1}) &\mbox{  } otherwise,\\
\end{cases}$\\
is a uninorm on $[b_{i},1]$, where $U_{1}^{T}=S_{e}$.

\end{theorem}

\proof
Let $(L,\leq,0,1)$ be a bounded lattice with   $a_{0},a_{1},a_{2},\ldots,a_{n}\in L$ such that $0_{L}=a_{0}\leq a_{1}\leq a_{2}\leq \ldots \leq a_{n}=1$.  Then we give the proof of the fact that $U_{i}^{S}$ is a uninorm. The same result for $U_{i}^{T}$ can be obtained using similar arguments.

If we take $i=2$  in the definition of  $U_{i}^{S}$,  then we can see that  $U_{2}^{S}$ is exactly  the uninorm $U_{d}$ in Theorem \ref{th2.5}.

Assume that $U_{k}^{S}$ is a uninorm. Then we prove that $U_{k+1}^{S}$ is a uninorm, i.e.

$U_{k+1}^{S}(x,y)=\begin{cases}
U_{k}^{S}(x, y) &\mbox{if } (x,y)\in [0,a_{k}]^{2},\\
S_{k+1}(x, y) &\mbox{if } (x,y)\in ]a_{k},a_{k+1}]^{2},\\
x\vee y &\mbox{if } (x,y)\in [0,e]\times (I_{a_{k}}^{e}\cup]a_{k},a_{k+1}])\\
&\mbox{ } \ \ \ \ \ \ \ \ \ \ \ \cup (I_{a_{k}}^{e}\cup]a_{k},a_{k+1}])\times[0,e],\\
x &\mbox{if } (x,y)\in I_{e,a_{k}}\times [0,e],\\
y &\mbox{if } (x,y)\in [0,e]\times I_{e,a_{k}},\\
S_{k+1}(x\vee a _{k}, y\vee a _{k}) &\mbox{  } otherwise.\\
\end{cases}$\\

The commutativity of $U_{k+1}^{S}$  and the fact that $e$  is a neutral element are clearly ture. Hence, we prove only the monotonicity and the associativity of $U_{k+1}^{S}$.

I. Monotonicity: We prove that if $x\leq y$, then $U_{k+1}^{S}(x,z)\leq U_{k+1}^{S}(y,z)$ for all $z\in L$.  The proof is split into all possible cases:

1. Let $x\in [0,e]$.

\ \ \ 1.1. $y\in [0,e]$,

\ \ \ \ \ \ 1.1.1. $z\in [0,e]\cup I_{e}^{a_{k}}\cup ]e,a_{k}]$,

\ \ \ \ \ \ \ \ \ \ \ \ Since the monotonicity property of $U_{k}^{S}(x,y)$, the monotonicity is true.

\ \ \ \ \ \ 1.1.2. $z\in I_{a_{k}}^{e}\cup I_{e,a_{k}}\cup ]a_{k},a_{k+1}]$,

\ \ \ \ \ \ \ \ \ \ \ \ $U_{k+1}^{S}(x,z)=z =U_{k+1}^{S}(y,z)$

\ \ \ 1.2. $y\in I_{e}^{a_{k}}$,

\ \ \ \ \ \ 1.2.1. $z\in [0,e]\cup I_{e}^{a_{k}}\cup ]e,a_{k}]$,

\ \ \ \ \ \ \ \ \ \ \ \ Since the monotonicity property of $U_{k}^{S}(x,y)$, the monotonicity is true.

\ \ \ \ \ \ 1.2.2. $z\in I_{a_{k}}^{e}\cup I_{e,a_{k}}$,

\ \ \ \ \ \ \ \ \ \ \ \ $U_{k+1}^{S}(x,z)=z<z\vee a_{k}=S_{k+1}(y\vee a_{k},z\vee a_{k})=U_{k+1}^{S}(y,z)$

\ \ \ \ \ \ 1.2.3. $z\in ]a_{k},a_{k+1}]$,

\ \ \ \ \ \ \ \ \ \ \ \ $U_{k+1}^{S}(x,z)=z=S_{k+1}(y\vee a_{k},z\vee a_{k})=U_{k+1}^{S}(y,z)$

\ \ \ 1.3. $y\in ]e,a_{k}]$,

\ \ \ \ \ \ 1.3.1. $z\in [0,e]\cup I_{e}^{a_{k}}\cup ]e,a_{k}]$,

\ \ \ \ \ \ \ \ \ \ \ \ Since the monotonicity property of $U_{k}^{S}(x,y)$, the monotonicity is true.

\ \ \ \ \ \ 1.3.2. $z\in I_{a_{k}}^{e}\cup I_{e,a_{k}}$,

\ \ \ \ \ \ \ \ \ \ \ \ $U_{k+1}^{S}(x,z)=z<z\vee a_{k}=S_{k+1}(y\vee a_{k},z\vee a_{k})=U_{k+1}^{S}(y,z)$

\ \ \ \ \ \ 1.3.3. $z\in ]a_{k},a_{k+1}]$,

\ \ \ \ \ \ \ \ \ \ \ \ $U_{k+1}^{S}(x,z)=z=S_{k+1}(y\vee a_{k},z\vee a_{k})=U_{k+1}^{S}(y,z)$

\ \ \ 1.4. $y\in I_{a_{k}}^{e}\cup I_{e,a_{k}}$,

\ \ \ \ \ \ 1.4.1. $z\in [0,e]$,

\ \ \ \ \ \ \ \ \ \ \ \ $U_{k+1}^{S}(x,z)=U_{k}^{S}(x,z)\leq x<y=U_{k+1}^{S}(y,z)$

\ \ \ \ \ \ 1.4.2. $z\in I_{e}^{a_{k}}\cup ]e,a_{k}]$,

\ \ \ \ \ \ \ \ \ \ \ \ $U_{k+1}^{S}(x,z)=U_{k}^{S}(x,z)\leq a_{k}<y\vee a_{k}=S_{k+1}(y\vee a_{k},z\vee a_{k})=U_{k+1}^{S}(y,z)$

\ \ \ \ \ \ 1.4.3. $z\in I_{a_{k}}^{e}\cup I_{e,a_{k}}$,

\ \ \ \ \ \ \ \ \ \ \ \ $U_{k+1}^{S}(x,z)=z< z\vee a_{k}\leq S_{k+1}(y\vee a_{k},z\vee a_{k})=U_{k+1}^{S}(y,z)$

\ \ \ \ \ \ 1.4.4. $z\in ]a_{k},a_{k+1}]$,

\ \ \ \ \ \ \ \ \ \ \ \ $U_{k+1}^{S}(x,z)=z\leq S_{k+1}(y\vee a_{k},z)=U_{k+1}^{S}(y,z)$

\ \ \ 1.5. $y\in ]a_{k},a_{k+1}]$,

\ \ \ \ \ \ 1.5.1. $z\in [0,e]$,

\ \ \ \ \ \ \ \ \ \ \ \ $U_{k+1}^{S}(x,z)=U_{k}^{S}(x,z)<y=U_{k+1}^{S}(y,z)$

\ \ \ \ \ \ 1.5.2. $z\in  I_{e}^{a_{k}}\cup ]e,a_{k}]$,

\ \ \ \ \ \ \ \ \ \ \ \ $U_{k+1}^{S}(x,z)=U_{K}^{S}(x,z)<y=S_{k+1}(y\vee a_{k},z\vee a_{k})=U_{k+1}^{S}(y,z)$

\ \ \ \ \ \ 1.5.3. $z\in I_{a_{k}}^{e}\cup I_{e,a_{k}}$,

\ \ \ \ \ \ \ \ \ \ \ \ $U_{k+1}^{S}(x,z)=z<z\vee a_{k}\leq S_{k+1}(y,z\vee a_{k})=U_{k+1}^{S}(y,z)$

\ \ \ \ \ \ 1.5.4. $z\in ]a_{k},a_{k+1}]$,

\ \ \ \ \ \ \ \ \ \ \ \ $U_{k+1}^{S}(x,z)=z\leq S_{k+1}(y,z)=U_{k+1}^{S}(y,z)$

2. Let $x\in I_{e}^{a_{k}}$.

\ \ \ 2.1. $y\in I_{e}^{a_{k}}$,

\ \ \ \ \ \ 2.1.1. $z\in [0,e]\cup I_{e}^{a_{k}}\cup ]e,a_{k}]$,

\ \ \ \ \ \ \ \ \ \ \ \ Since the monotonicity property of $U_{k}^{S}(x,y)$, the monotonicity is true.

\ \ \ \ \ \ 2.1.2. $z\in I_{a_{k}}^{e}\cup I_{e,a_{k}}$,

\ \ \ \ \ \ \ \ \ \ \ \ $U_{k+1}^{S}(x,z)=S_{k+1}(x\vee a_{k},z\vee a_{k})=z\vee a_{k}=S_{k+1}(y\vee a_{k},z\vee a_{k})=U_{k+1}^{S}(y,z)$

\ \ \ \ \ \ 2.1.3. $z\in ]a_{k},a_{k+1}]$,

\ \ \ \ \ \ \ \ \ \ \ \ $U_{k+1}^{S}(x,z)=S_{k+1}(x\vee a_{k},z\vee a_{k})=z=S_{k+1}(y\vee a_{k},z\vee a_{k})=U_{k+1}^{S}(y,z)$

\ \ \ 2.2. $y\in ]e,a_{k}]$,

\ \ \ \ \ \ 2.2.1. $z\in [0,e]\cup I_{e}^{a_{k}}\cup ]e,a_{k}]$,

\ \ \ \ \ \ \ \ \ \ \ \ Since the monotonicity property of $U_{k}^{S}(x,y)$, the monotonicity is true.

\ \ \ \ \ \ 2.2.2. $z\in I_{a_{k}}^{e}\cup I_{e,a_{k}}$,

\ \ \ \ \ \ \ \ \ \ \ \ $U_{k+1}^{S}(x,z)=S_{k+1}(x\vee a_{k},z\vee a_{k})=z\vee a_{k}=S_{k+1}(y\vee a_{k},z\vee a_{k})=U_{k+1}^{S}(y,z)$

\ \ \ \ \ \ 2.2.3. $z\in ]a_{k},a_{k+1}]$,

\ \ \ \ \ \ \ \ \ \ \ \ $U_{k+1}^{S}(x,z)=S_{k+1}(x\vee a_{k},z\vee a_{k})=z=S_{k+1}(y\vee a_{k},z\vee a_{k})=U_{k+1}^{S}(y,z)$

\ \ \ 2.3. $y\in I_{a_{k}}^{e}\cup I_{e,a_{k}}$,

\ \ \ \ \ \ 2.3.1. $z\in [0,e]$,

\ \ \ \ \ \ \ \ \ \ \ \ $U_{k+1}^{S}(x,z)=U_{k}^{S}(x,z)\leq x<y=U_{k+1}^{S}(y,z)$

\ \ \ \ \ \ 2.3.2. $z\in I_{e}^{a_{k}}\cup ]e,a_{k}]$,

\ \ \ \ \ \ \ \ \ \ \ \ $U_{k+1}^{S}(x,z)=U_{k}^{S}(x,z)\leq a_{k}<y\vee a_{k}=S_{k+1}(y\vee a_{k},z\vee a_{k})=U_{k+1}^{S}(y,z)$

\ \ \ \ \ \ 2.3.3. $z\in I_{a_{k}}^{e}\cup I_{e,a_{k}}$,

\ \ \ \ \ \ \ \ \ \ \ \ $U_{k+1}^{S}(x,z)=S_{k+1}(x\vee a_{k},z\vee a_{k})=z\vee a_{k}\leq S_{k+1}(y\vee a_{k},z\vee a_{k})=U_{k+1}^{S}(y,z)$

\ \ \ \ \ \ 2.3.4. $z\in ]a_{k},a_{k+1}]$,

\ \ \ \ \ \ \ \ \ \ \ \ $U_{k+1}^{S}(x,z)=z\leq S_{k+1}(y\vee a_{k},z)=U_{k+1}^{S}(y,z)$

\ \ \ 2.4. $y\in ]a_{k},a_{k+1}]$,

\ \ \ \ \ \ 2.4.1. $z\in [0,e]$,

\ \ \ \ \ \ \ \ \ \ \ \ $U_{k+1}^{S}(x,z)=U_{k}^{S}(x,z)<y=U_{k+1}^{S}(y,z)$

\ \ \ \ \ \ 2.4.2. $z\in  I_{e}^{a_{k}}\cup ]e,a_{k}]$,

\ \ \ \ \ \ \ \ \ \ \ \ $U_{k+1}^{S}(x,z)=U_{k}^{S}(x,z)<y=S_{k+1}(y\vee a_{k},z\vee a_{k})=U_{k+1}^{S}(y,z)$

\ \ \ \ \ \ 2.4.3. $z\in I_{a_{k}}^{e}\cup I_{e,a_{k}}$,

\ \ \ \ \ \ \ \ \ \ \ \ $U_{k+1}^{S}(x,z)=S_{k+1}(x\vee a_{k},z\vee a_{k})=z\vee a_{k}\leq S_{k+1}(y\vee a_{k},z\vee a_{k})=U_{k+1}^{S}(y,z)$

\ \ \ \ \ \ 2.4.4. $z\in ]a_{k},a_{k+1}]$,

\ \ \ \ \ \ \ \ \ \ \ \ $U_{k+1}^{S}(x,z)=S_{k+1}(x\vee a_{k},z\vee a_{k})=z\leq S_{k+1}(y,z)=U_{k+1}^{S}(y,z)$

3. Let $x\in ]e,a_{k}]$.

\ \ \ 3.1. $y\in ]e,a_{k}]$,

\ \ \ \ \ \ 3.1.1. $z\in [0,e]\cup I_{e}^{a_{k}}\cup ]e,a_{k}]$,

\ \ \ \ \ \ \ \ \ \ \ \ Since the monotonicity property of $U_{k}^{S}(x,y)$, the monotonicity is true.

\ \ \ \ \ \ 3.1.2. $z\in I_{a_{k}}^{e}\cup I_{e,a_{k}}$,

\ \ \ \ \ \ \ \ \ \ \ \ $U_{k+1}^{S}(x,z)=S_{k+1}(x\vee a_{k},z\vee a_{k})=z\vee a_{k}=S_{k+1}(y\vee a_{k},z\vee a_{k})=U_{k+1}^{S}(y,z)$

\ \ \ \ \ \ 3.1.3. $z\in ]a_{k},a_{k+1}]$,

\ \ \ \ \ \ \ \ \ \ \ \ $U_{k+1}^{S}(x,z)=S_{k+1}(x\vee a_{k},z\vee a_{k})=z=S_{k+1}(y\vee a_{k},z\vee a_{k})=U_{k+1}^{S}(y,z)$

\ \ \ 3.2. $y\in I_{a_{k}}^{e}$,

\ \ \ \ \ \ 3.2.1. $z\in [0,e]$,

\ \ \ \ \ \ \ \ \ \ \ \ $U_{k+1}^{S}(x,z)=U_{k}^{S}(x,z)\leq x<y=U_{k+1}^{S}(y,z)$

\ \ \ \ \ \ 3.2.2. $z\in I_{e}^{a_{k}}\cup ]e,a_{k}]$,

\ \ \ \ \ \ \ \ \ \ \ \ $U_{k+1}^{S}(x,z)=U_{k}^{S}(x,z)\leq a_{k}<y\vee a_{k}=S_{k+1}(y\vee a_{k},z\vee a_{k})=U_{k+1}^{S}(y,z)$

\ \ \ \ \ \ 3.2.3. $z\in I_{a_{k}}^{e}\cup I_{e,a_{k}}$,

\ \ \ \ \ \ \ \ \ \ \ \ $U_{k+1}^{S}(x,z)=S_{k+1}(x\vee a_{k},z\vee a_{k})=z\vee a_{k}\leq S_{k+1}(y\vee a_{k},z\vee a_{k})=U_{k+1}^{S}(y,z)$

\ \ \ \ \ \ 3.2.4. $z\in ]a_{k},a_{k+1}]$,

\ \ \ \ \ \ \ \ \ \ \ \ $U_{k+1}^{S}(x,z)=S_{k+1}(x\vee a_{k},z\vee a_{k})=z\leq S_{k+1}(y\vee a_{k},z)=U_{k+1}^{S}(y,z)$

\ \ \ 3.3. $y\in ]a_{k},a_{k+1}]$,

\ \ \ \ \ \ 3.3.1. $z\in [0,e]$,

\ \ \ \ \ \ \ \ \ \ \ \ $U_{k+1}^{S}(x,z)=U_{k}^{S}(x,z)\leq a_{k}<y=U_{k+1}^{S}(y,z)$

\ \ \ \ \ \ 3.3.2. $z\in  I_{e}^{a_{k}}\cup ]e,a_{k}]$,

\ \ \ \ \ \ \ \ \ \ \ \ $U_{k+1}^{S}(x,z)=U_{k}^{S}(x,z)\leq a_{k}<y=S_{k+1}(y\vee a_{k},z\vee a_{k})=U_{k+1}^{S}(y,z)$

\ \ \ \ \ \ 3.3.3. $z\in I_{a_{k}}^{e}\cup I_{e,a_{k}}$,

\ \ \ \ \ \ \ \ \ \ \ \ $U_{k+1}^{S}(x,z)=S_{k+1}(x\vee a_{k},z\vee a_{k})=z\vee a_{k}\leq S_{k+1}(y\vee a_{k},z\vee a_{k})=U_{k+1}^{S}(y,z)$

\ \ \ \ \ \ 3.4.4. $z\in ]a_{k},a_{k+1}]$,

\ \ \ \ \ \ \ \ \ \ \ \ $U_{k+1}^{S}(x,z)=S_{k+1}(x\vee a_{k},z\vee a_{k})=z\leq S_{k+1}(y,z)=U_{k+1}^{S}(y,z)$

4. Let $x\in I_{a_{k}}^{e}$.

\ \ \ 4.1. $y\in I_{a_{k}}^{e}$,

\ \ \ \ \ \ 4.1.1. $z\in [0,e]$,

\ \ \ \ \ \ \ \ \ \ \ \ $U_{k+1}^{S}(x,z)=x\leq y=U_{k+1}^{S}(y,z)$

\ \ \ \ \ \ 4.1.2. $z\in I_{e}^{a_{k}}\cup]e,a_{k}]\cup I_{a_{k}}^{e}\cup I_{e,a_{k}}\cup]a_{k},a_{k+1}]$,

\ \ \ \ \ \ \ \ \ \ \ \ Since the monotonicity property of $S_{k+1}$, the monotonicity is true.

\ \ \ 4.2. $y\in ]a_{k},a_{k+1}]$,

\ \ \ \ \ \ 4.2.1. $z\in [0,e]$,

\ \ \ \ \ \ \ \ \ \ \ \ $U_{k+1}^{S}(x,z)=x< y=U_{k+1}^{S}(y,z)$

\ \ \ \ \ \ 4.2.2. $z\in I_{e}^{a_{k}}\cup ]e,a_{k}]$,

\ \ \ \ \ \ \ \ \ \ \ \ $U_{k+1}^{S}(x,z)=S_{k+1}(x\vee a_{k},z\vee a_{k})=x\vee a_{k}\leq y=S_{k+1}(y\vee a_{k},z\vee a_{k})$

\ \ \ \ \ \ \ \ \ \ \ \ $=U_{k+1}^{S}(y,z)$

\ \ \ \ \ \ 4.2.3. $z\in I_{a_{k}}^{e}\cup I_{e,a_{k}}$,

\ \ \ \ \ \ \ \ \ \ \ \ $U_{k+1}^{S}(x,z)=S_{k+1}(x\vee a_{k},z\vee a_{k})\leq S_{k+1}(y,z\vee a_{k})=U_{k+1}^{S}(y,z)$

\ \ \ \ \ \ 4.2.4. $z\in ]a_{k},a_{k+1}]$,

\ \ \ \ \ \ \ \ \ \ \ \ $U_{k+1}^{S}(x,z)=S_{k+1}(x\vee a_{k},z\vee a_{k})=S_{k+1}(x\vee a_{k},z)\leq S_{k+1}(y,z)=U_{k+1}^{S}(y,z)$

5. Let $x\in I_{e,a_{k}}$.

\ \ \ 5.1. $y\in I_{e,a_{k}}$,

\ \ \ \ \ \ 5.1.1. $z\in [0,e]$,

\ \ \ \ \ \ \ \ \ \ \ \ $U_{k+1}^{S}(x,z)=x\leq y=U_{k+1}^{S}(y,z)$

\ \ \ \ \ \ 5.1.2. $z\in I_{e}^{a_{k}}\cup]e,a_{k}]\cup I_{a_{k}}^{e}\cup I_{e,a_{k}}\cup]a_{k},a_{k+1}]$,

\ \ \ \ \ \ \ \ \ \ \ \ Since the monotonicity property of $S_{k+1}$, the monotonicity is true.

\ \ \ 5.2. $y\in I_{a_{k}}^{e}$,

\ \ \ \ \ \ 5.2.1. $z\in [0,e]$,

\ \ \ \ \ \ \ \ \ \ \ \ $U_{k+1}^{S}(x,z)=x< y=U_{k+1}^{S}(y,z)$

\ \ \ \ \ \ 5.2.2. $z\in I_{e}^{a_{k}}\cup]e,a_{k}]\cup I_{a_{k}}^{e}\cup I_{e,a_{k}}\cup]a_{k},a_{k+1}]$,

\ \ \ \ \ \ \ \ \ \ \ \ Since the monotonicity property of $S_{k+1}$, the monotonicity is true.

\ \ \ 5.3. $y\in ]a_{k},a_{k+1}]$,

\ \ \ \ \ \ 5.3.1. $z\in [0,e]$,

\ \ \ \ \ \ \ \ \ \ \ \ $U_{k+1}^{S}(x,z)=x<y=U_{k+1}^{S}(y,z)$

\ \ \ \ \ \ 5.3.2. $z\in I_{e}^{a_{k}}\cup ]e,a_{k}]$,

\ \ \ \ \ \ \ \ \ \ \ \ $U_{k+1}^{S}(x,z)=S_{k+1}(x\vee a_{k},z\vee a_{k})=x\vee a_{k}\leq y=S_{k+1}(y\vee a_{k},z\vee a_{k})=U_{k+1}^{S}(y,z)$

\ \ \ \ \ \ 5.3.3. $z\in I_{a_{k}}^{e}\cup I_{e,a_{k}}$,

\ \ \ \ \ \ \ \ \ \ \ \ $U_{k+1}^{S}(x,z)=S_{k+1}(x\vee a_{k},z\vee a_{k})\leq S_{k+1}(y\vee a_{k},z\vee a_{k})=U_{k+1}^{S}(y,z)$

\ \ \ \ \ \ 5.3.4. $z\in ]a_{k},a_{k+1}]$,

\ \ \ \ \ \ \ \ \ \ \ \ $U_{k+1}^{S}(x,z)=S_{k+1}(x\vee a_{k},z\vee a_{k})=S_{k+1}(x\vee a_{k},z)\leq S_{k+1}(y,z)=U_{k+1}^{S}(y,z)$

6. Let $x\in ]a_{k},a_{k+1}]$.

\ \ \ 6.1. $y\in ]a_{k},a_{k+1}]$,

\ \ \ \ \ \ 6.1.1. $z\in [0,e]$,

\ \ \ \ \ \ \ \ \ \ \ \ $U_{k+1}^{S}(x,z)=x\leq y=U_{k+1}^{S}(y,z)$

\ \ \ \ \ \ 6.1.2. $z\in I_{e}^{a_{k}}\cup]e,a_{k}]\cup I_{a_{k}}^{e}\cup I_{e,a_{k}}\cup]a_{k},a_{k+1}]$,

\ \ \ \ \ \ \ \ \ \ \ \ Since the monotonicity property of $S_{k+1}$, the monotonicity is true.

Combining the above cases, we obtain that $U_{k+1}^{S}(x,z)\leq U_{k+1}^{S}(y,z)$ holds for all $x,y,z\in L$ such that $x\leq y$. Therefore, $U_{k+1}^{S}$ is monotonic.

II. Associativity: It can be shown that $U_{k+1}^{S}(x,U_{k+1}^{S}(y,z))=U_{k+1}^{S}(U_{k+1}^{S}(x,y),z)$  for all $x,y,z\in L$. By Proposition  \ref{pro2.2}, we only need to consider the following cases:

1. If $x,y,z\in [0,e]\cup I_{e}^{a_{k}}\cup]e,a_{k}]  $, then since $U_{k}^{S}$ is associative, $U_{k+1}^{S}$ is associative.

2. If $x,y,z\in I_{a_{k}}^{e}\cup I_{e,a_{k}}\cup ]a_{k},a_{k+1}] $, then since $S_{k+1}$ is associative, $U_{k+1}^{S}$ is associative.

3. If $x,y\in [0,e],z\in I_{e}^{a_{k}}\cup  ]e,a_{k}] $, then since $U_{k}^{S}$ is associative, $U_{k+1}^{S}$ is associative.

4. If $x,y\in [0,e],z\in I_{a_{k}}^{e}\cup I_{e,a_{k}}\cup ]a_{k},a_{k+1}] $, then
$U_{k+1}^{S}(x,U_{k+1}^{S}(y,z))
=U_{k+1}^{S}(x,z)=z=U_{k+1}^{S}(U_{k}^{S}(x,y),z)=U_{k+1}^{S}(U_{k+1}^{S}(x,y),z).$ 
Thus $U_{k+1}^{S}$ is associative.


5. If $x,y\in I_{e}^{a_{k}},z\in ]e,a_{k}]$, then since $U_{k}^{S}$ is associative, $U_{k+1}^{S}$ is associative.


6. If $x,y\in I_{e}^{a_{k}},z\in I_{a_{k}}^{e}\cup I_{e,a_{k}} $, then  $U_{k+1}^{S}(x,U_{k+1}^{S}(y,z))=U_{k+1}^{S}(x,S_{k+1}(y\vee a_{k},z\vee a_{k}))=U_{k+1}^{S}(x,z\vee a_{k})=z\vee a_{k}=S_{k+1}(U_{k}^{S}(x,y)\vee a_{k},z\vee a_{k})=U_{k+1}^{S}(U_{k}^{S}(x,y),z)=U_{k+1}^{S}(U_{k+1}^{S}(x,y),z).$
Thus $U_{k+1}^{S}$ is associative.

7. If $x,y\in I_{e}^{a_{k}},z\in ]a_{k},a_{k+1}]$, then $U_{k+1}^{S}(x,U_{k+1}^{S}(y,z))=U_{k+1}^{S}(x,S_{k+1}(y\vee a_{k},z\vee a_{k}))=U_{k+1}^{S}(x,z)=S_{k+1}(x\vee a_{k},z\vee a_{k})=z$ and $U_{k+1}^{S}(U_{k+1}^{S}(x,y),z)=U_{k+1}^{S}(U_{k}^{S}(x,y),z)=S_{k+1}(U_{k}^{S}(x,y)\vee a_{k},z\vee a_{k})=z$. Thus $U_{k+1}^{S}$ is associative.

8. If $x,y\in ]e,a_{k}],z\in I_{a_{k}}^{e} \cup  I_{e,a_{k}}$, then $U_{k+1}^{S}(x,U_{k+1}^{S}(y,z))=U_{k+1}^{S}(x,S_{k+1}(y\vee a_{k},z\vee a_{k}))=U_{k+1}^{S}(x,z\vee a_{k})=z\vee a_{k}$ and $U_{k+1}^{S}(U_{k+1}^{S}(x,y),z)=U_{k+1}^{S}(U_{k}^{S}(x,y),z)=S_{k+1}(U_{k}^{S}(x,y)\\\vee a_{k},z\vee a_{k})=z\vee a_{k}$. Thus $U_{k+1}^{S}$ is associative.

9. If $x,y\in ]e,a_{k}],z\in ]a_{k},a_{k+1}]$, then $U_{k+1}^{S}(x,U_{k+1}^{S}(y,z))=U_{k+1}^{S}(x,S_{k+1}(y\vee a_{k},z\vee a_{k}))=U_{k+1}^{S}(x,z)=S_{k+1}(x\vee a_{k},z\vee a_{k})=z$ and $U_{k+1}^{S}(U_{k+1}^{S}(x,y),z)=U_{k+1}^{S}(U_{k}^{S}(x,y),z)=S_{k+1}(U_{k}^{S}(x,y)\vee a_{k},z\vee a_{k})=z$. Thus $U_{k+1}^{S}$ is associative.

10. If $x,y\in I_{a_{k}}^{e},z\in I_{e,a_{k}}\cup ]a_{k},a_{k+1}]$, then since $S_{k+1}$ is associative,
$U_{k+1}^{S}$ is associative.

11. If $x,y\in I_{e,a_{k}},z\in ]a_{k},a_{k+1}]$,  then since $S_{k+1}$ is associative, we have  $U_{k+1}^{S}$ is associative.

12. If $x\in [0,e], y,z\in I_{e}^{a_{k}}\cup ]e,a_{k}] $, then since $U_{k}^{S}$ is associative, $U_{k+1}^{S}$ is associative.


13. If $x\in [0,e],y,z\in I_{a_{k}}^{e}\cup I_{e,a_{k}} $, then $U_{k+1}^{S}(x,U_{k+1}^{S}(y,z))=U_{k+1}^{S}(x,S_{k+1}(y\vee a_{k},z\vee a_{k}))=S_{k+1}(y\vee a_{k},z\vee a_{k})=U_{k+1}^{S}(y,z)=U_{k+1}^{S}(U_{k+1}^{S}(x,y),z).$  
Thus $U_{k+1}^{S}$ is associative.

14. If $x\in [0,e],y,z\in ]a_{k},a_{k+1}]$, then $U_{k+1}^{S}(x,U_{k+1}^{S}(y,z))=U_{k+1}^{S}(x,S_{k+1}(y,z))
=S_{k+1}(y,z)=U_{k+1}^{S}(y,z)=U_{k+1}^{S}(U_{k+1}^{S}(x,y),z).$  
Thus $U_{k+1}^{S}$ is associative.

15. If $x\in I_{e}^{a_{k}},y,z\in ]e,a_{k}]$, then since $U_{k}^{S}$ is associative, $U_{k+1}^{S}$ is associative.

16. If $x\in I_{e}^{a_{k}},y,z\in I_{a_{k}}^{e}\cup I_{e,a_{k}}\cup ]a_{k},a_{k+1}] $,  then since $S_{k+1}$ is associative,   $U_{k+1}^{S}$ is associative.

17. If $x\in ]e,a_{k}],y,z\in I_{a_{k}}^{e}\cup I_{e,a_{k}}\cup ]a_{k},a_{k+1}] $,  then since $S_{k+1}$ is associative,    $U_{k+1}^{S}$ is associative.

18. If $x\in I_{a_{k}}^{e},y,z\in I_{e,a_{k}}\cup ]a_{k},a_{k+1}]$, then since $S_{k+1}$ is associative,   $U_{k+1}^{S}$ is associative.

19. If $x\in I_{e,a_{k}},y,z\in ]a_{k},a_{k+1}]$, then then since $S_{k+1}$ is associative,    $U_{k+1}^{S}$ is associative.

20. If $x\in [0,e],y\in I_{e}^{a_{k}}, z\in ]e,a_{k}]$, then since $U_{k}^{S}$ is associative, $U_{k+1}^{S}$ is associative.

21. If $x\in [0,e],y\in I_{e}^{a_{k}},z\in I_{a_{k}}^{e}\cup I_{e,a_{k}} $, then $U_{k+1}^{S}(x,U_{k+1}^{S}(y,z))=U_{k+1}^{S}(x,S_{k+1}(a_{k},z\vee a_{k}))=z\vee a_{k}=S_{k+1}(a_{k},z\vee a_{k})=S_{k+1}(U_{k}^{S}(x,y)\vee a_{k},z\vee a_{k})=U_{k+1}^{S}(U_{k}^{S}(x,y),z)=U_{k+1}^{S}(U_{k+1}^{S}(x,y),z)$ 
and $U_{k+1}^{S}(y,U_{k+1}^{S}(x,z))=U_{k+1}^{S}(y,z)=S_{k+1}(a_{k},z\vee a_{k})=z\vee a_{k}$. Thus $U_{k+1}^{S}$ is associative.

22. If $x\in [0,e],y\in I_{e}^{a_{k}},z\in ]a_{k},a_{k+1}]$, then $U_{k+1}^{S}(x,U_{k+1}^{S}(y,z))=U_{k+1}^{S}(x,S_{k+1}(y\vee a_{k},z\vee a_{k}))=U_{k+1}^{S}(x,z)=z,U_{k+1}^{S}(U_{k+1}^{S}(x,y),z)=U_{k+1}^{S}(U_{k}^{S}(x,y),z)=S_{k+1}(U_{k}^{S}(x,y)\vee a_{k},z\vee a_{k})=z$ and $U_{k+1}^{S}(y,U_{k+1}^{S}(x,z))=U_{k+1}^{S}(y,z)=S_{k+1}(y\vee a_{k},z\vee a_{k}) =z$. Thus $U_{k+1}^{S}$ is associative.

23. If $x\in [0,e],y\in ]e,a_{k}],z\in I_{a_{k}}^{e}\cup I_{e,a_{k}} $, then $U_{k+1}^{S}(x,U_{k+1}^{S}(y,z))=U_{k+1}^{S}(x,S_{k+1}(y\vee a_{k},z\vee a_{k}))=S_{k+1}(y\vee a_{k},z\vee a_{k})=z\vee a_{k},U_{k+1}^{S}(U_{k+1}^{S}(x,y),z)=U_{k+1}^{S}(U_{k}^{S}(x,y),z)=S_{k+1}(U_{k}^{S}(x,y)\vee a_{k},z\vee a_{k})=z\vee a_{k}$ and $U_{k+1}^{S}(y,U_{k+1}^{S}(x,z))=U_{k+1}^{S}(y,z)=S_{k+1}(y\vee a_{k},z\vee a_{k})=z\vee a_{k}$. Thus $U_{k+1}^{S}$ is associative.

24. If $x\in [0,e],y\in ]e,a_{k}],z\in ]a_{k},a_{k+1}]$, then $U_{k+1}^{S}(x,U_{k+1}^{S}(y,z))=U_{k+1}^{S}(x,S_{k+1}(y\vee a_{k},z\vee a_{k}))=U_{k+1}^{S}(x,z)=z,U_{k+1}^{S}(U_{k+1}^{S}(x,y),z)=U_{k+1}^{S}(U_{k}^{S}(x,y),z)=S_{k+1}(U_{k}^{S}(x,y)\vee a_{k},z\vee a_{k})=z$ and $U_{k+1}^{S}(y,U_{k+1}^{S}(x,z))=U_{k+1}^{S}(y,z)=S_{k+1}(y\vee a_{k},z\vee a_{k}) =z$. Thus $U_{k+1}^{S}$ is associative.

25. If $x\in [0,e],y\in I_{a_{k}}^{e},z\in I_{e,a_{k}}$, then
 $U_{k+1}^{S}(x,U_{k+1}^{S}(y,z))=U_{k+1}^{S}(x,S_{k+1}(y\vee a_{k},z\vee a_{k}))=S_{k+1}(y\vee a_{k},z\vee a_{k}),U_{k+1}^{S}(U_{k+1}^{S}(x,y),z)=U_{k+1}^{S}(y,z)=S_{k+1}(y\vee a_{k},z\vee a_{k})$ and $U_{k+1}^{S}(y,U_{k+1}^{S}(x,z))=U_{k+1}^{S}(y,z)=S_{k+1}(y\vee a_{k},z\vee a_{k})$. Thus $U_{k+1}^{S}$ is associative.

26. If $x\in [0,e],y\in I_{a_{k}}^{e},z\in ]a_{k},a_{k+1}]$, then $U_{k+1}^{S}(x,U_{k+1}^{S}(y,z))=U_{k+1}^{S}(x,S_{k+1}(y\vee a_{k},z\vee a_{k}))=S_{k+1}(y\vee a_{k},z),U_{k+1}^{S}(U_{k+1}^{S}(x,y),z)=U_{k+1}^{S}(y,z)=S_{k+1}(y\vee a_{k},z)$ and $U_{k+1}^{S}(y,U_{k+1}^{S}(x,z))=U_{k+1}^{S}(y,z)=S_{k+1}(y\vee a_{k},z)$. Thus $U_{k+1}^{S}$ is associative.

27. If $x\in [0,e],y\in I_{e,a_{k}},z\in ]a_{k},a_{k+1}]$, then $U_{k+1}^{S}(x,U_{k+1}^{S}(y,z))=U_{k+1}^{S}(x,S_{k+1}(y\vee a_{k},z\vee a_{k}))=S_{k+1}(y\vee a_{k},z),U_{k+1}^{S}(U_{k+1}^{S}(x,y),z)=U_{k+1}^{S}(y,z)=S_{k+1}(y\vee a_{k},z)$ and $U_{k+1}^{S}(y,U_{k+1}^{S}(x,z))=U_{k+1}^{S}(y,z)=S_{k+1}(y\vee a_{k},z)$. Thus $U_{k+1}^{S}$ is associative.

28. If $x\in I_{e}^{a_{k}},y\in ]e,a_{k}],z\in I_{a_{k}}^{e}\cup I_{e,a_{k}} $, then $U_{k+1}^{S}(x,U_{k+1}^{S}(y,z))=U_{k+1}^{S}(x,S_{k+1}(y\vee a_{k},z\vee a_{k}))=U_{k+1}^{S}(x,z\vee a_{k})=z\vee a_{k}=S_{k+1}(U_{k}^{S}(x,y)\vee a_{k},z\vee a_{k})=U_{k+1}^{S}(U_{k}^{S}(x,y),z)=U_{k+1}^{S}(U_{k+1}^{S}(x,y),z)$
and $U_{k+1}^{S}(y,U_{k+1}^{S}(x,z))=U_{k+1}^{S}(y,S_{k+1}(x\vee a_{k},z\vee a_{k}))=S_{k+1}(x\vee a_{k},z\vee a_{k})=z\vee a_{k}$. Thus $U_{k+1}^{S}$ is associative.

29. If $x\in I_{e}^{a_{k}},y\in ]e,a_{k}],z\in ]a_{k},a_{k+1}]$, then $U_{k+1}^{S}(x,U_{k+1}^{S}(y,z))=U_{k+1}^{S}(x,S_{k+1}(y\vee a_{k},z\vee a_{k}))=U_{k+1}^{S}(x,z\vee a_{k})=S_{k+1}(x\vee a_{k},z\vee a_{k})=z=S_{k+1}(U_{k}^{S}(x,y)\vee a_{k},z\vee a_{k})=U_{k+1}^{S}(U_{k}^{S}(x,y),z)=U_{k+1}^{S}(U_{k+1}^{S}(x,y),z)$
and $U_{k+1}^{S}(y,U_{k+1}^{S}(x,z))=U_{k+1}^{S}(y,S_{k+1}(x\vee a_{k},z\vee a_{k}))=U_{k+1}^{S}(y,z)=S_{k+1}(y\vee a_{k},z\vee a_{k})=z$. Thus $U_{k+1}^{S}$ is associative.

30. If $x\in I_{e}^{a_{k}},y\in I_{a_{k}}^{e},z\in I_{e,a_{k}}\cup]a_{k},a_{k+1}]$,  then since $S_{k+1}$ is associative,   $U_{k+1}^{S}$ is associative.

31. If $x\in I_{e}^{a_{k}},y\in I_{e,a_{k}},z\in ]a_{k},a_{k+1}]$,  then since $S_{k+1}$ is associative,  $U_{k+1}^{S}$ is associative.

32. If $x\in ]e,a_{k}],y\in I_{a_{k}}^{e},z\in I_{e,a_{k}}\cup ]a_{k},a_{k+1}]$,  then since $S_{k+1}$ is associative, $U_{k+1}^{S}$ is associative.

33. If $x\in ]e,a_{k}],y\in I_{e,a_{k}},z\in ]a_{k},a_{k+1}]$,  then since $S_{k+1}$ is associative, $U_{k+1}^{S}$ is associative.

34. If $x\in I_{a_{k}}^{e},y\in I_{e,a_{k}},z\in ]a_{k},a_{k+1}]$,  then since $S_{k+1}$ is associative,    $U_{k+1}^{S}$ is associative.

Form $1$ to $34$, we obtain that $U_{k+1}^{S}(x,U_{k+1}^{S}(y,z))=U_{k+1}^{S}(U_{k+1}^{S}(x,y),z)$ for all $x,y,z\in L$ by Proposition \ref{pro2.2}. Therefore, $U_{k+1}^{S}$ is a uninorm on $L$ with the neutral element $e$.

Thus, $U_{i}^{S}$ is a uninorm on $L$ with the neutral element $e$ for $i\in \{2,3,\ldots,n\}$.

\begin{corollary}\label{co3.01}
Let $(L,\leq,0,1)$ be a bounded lattice.\\
(i) If we take $S_{i}=S_{\vee}$ on $[a_{i-1},a_{i}]$ of $L$  for $i\in \{2,3,\ldots,n\}$  in Theorem \ref{th3.0}, then we obtain the following uninorm on $L$.

$U_{i}^{\vee}(x,y)=\begin{cases}
U_{i-1}^{\vee}(x, y) &\mbox{if } (x,y)\in [0,a_{i-1}]^{2},\\

x\vee y &\mbox{if } (x,y)\in ]a_{i-1},a_{i}]^{2}\cup [0,e]\times (I_{a_{i-1}}^{e}\cup]a_{i-1},a_{i}])\\
&\mbox{ } \ \ \ \ \ \ \ \ \ \ \ \cup (I_{a_{i-1}}^{e}\cup]a_{i-1},a_{i}])\times[0,e],\\
x &\mbox{if } (x,y)\in I_{e,a_{i-1}}\times [0,e],\\
y &\mbox{if } (x,y)\in [0,e]\times I_{e,a_{i-1}},\\
x\vee y\vee a _{i-1} &\mbox{  } otherwise.\\
\end{cases}$

(ii) If we take $T_{i}=T_{\wedge}$ on $[b_{i},b_{i-1}]$ of $L$  for $i\in \{2,3,\ldots,n\}$ in Theorem \ref{th3.0}, then we obtain the following uninorm on $L$.

$U_{i}^{\wedge}(x,y)=\begin{cases}
U_{i-1}^{\wedge}(x, y) &\mbox{if } (x,y)\in [b_{i-1},1]^{2},\\
x\wedge y &\mbox{if } (x,y)\in [b_{i},b_{i-1})^{2}\cup[e,1]\times (I_{b_{i-1}}^{e}\cup[b_{i},b_{i-1}[)\\
&\mbox{ } \ \ \ \ \ \ \ \ \ \ \ \cup (I_{b_{i-1}}^{e}\cup[b_{i},b_{i-1}[)\times[e,1],\\
x &\mbox{if } (x,y)\in I_{e,b_{i-1}}\times [e,1],\\
y &\mbox{if } (x,y)\in [e,1]\times I_{e,b_{i-1}},\\
x\wedge y\wedge b_{i-1} &\mbox{  } otherwise.\\
\end{cases}$\\

\end{corollary}

\begin{corollary}
Let $(L,\leq,0,1)$ be a bounded lattice.\\
(i) If we take $i=2$   in Corollary \ref{co3.01}, then we obtain the uninorm $U_{(T,e)}$ in Theorem \ref{th2.0}.\\
(ii) If we take $i=2$  in Corollary \ref{co3.01}, then we obtain the uninorm $U_{(S,e)}$ in Theorem \ref{th2.0}.
\end{corollary}

\section{\bf Conclusion}

In this article,  we still investigate the constructions of uninorms on   bounded lattices.

On one hand,  we investigate the constructions of uninorms on some appropriate bounded lattices with $e\in L\setminus\{0,1\}$, based on some additional constraints.   In this progress, we mainly focus on the  rationality of these constraints and then show that some of these constraints are sufficient and necessary for  the construction of uninorms.
Meanwhile, we study the constraints of some   uninorms and find that some constraints are  also necessary.  In the future research, when  considering the new  construction of uninorms,  if the constraints are needed, we expect they are sufficient and necessary.

 On the other hand,   we consider  new methods of constructing of uninorms, that is,  the iterative method using $t$-norms and $t$-conorms.  In the future research,  we will give more structures of  uninorms by this  method.


%



\end{document}